\documentclass{amsart}
\usepackage{amssymb}

\newcommand{\Q}{{\mathbb Q}}
\newcommand{\Z}{{\mathbb Z}}
\newcommand{\bN}{{\mathbb N}}
\newcommand{\N}{\mathcal N}

\newcommand{\R}{{\mathbb R}}
\newcommand{\F}{{\mathbb F}}
\newcommand{\K}{{\mathbb K}}

\newcommand{\LL}{{\mathbb L}}

\newcommand{\OO}{{\mathcal O}}

\newcommand{\pp}{{\mathcal{P}}}

\newcommand{\gS}{{\mathfrak S}}

\newcommand{\Tc}{{\mathcal T}}
\newcommand{\al}{{\alpha}}

\newcommand{\La}{{\Lambda}}
\newcommand{\eps}{{\varepsilon}}
\newcommand{\M}{{\rm M}}
\newcommand{\h}{{\rm h}}

\newcommand{\cI}{{\mathcal I}}
\newcommand{\cM}{{\mathcal M}}
\newcommand{\cN}{{\mathcal N}}

\newcommand{\bC}{{\mathbb C}}

\newcommand{\kro}[2]{\left( \frac{#1}{#2} \right) }

\begin {document}
\newtheorem{thm}{Theorem}
\newtheorem{lem}{Lemma}[section]
\newtheorem{theorem}[lem]{Theorem}

\newtheorem{prop}[lem]{Proposition}

\newtheorem*{remark}{Remark}
\theoremstyle{definition}

\theoremstyle{remark}

\title[Classical and Modular Approaches]{Classical
and Modular Approaches to Exponential Diophantine Equations \\
I. Fibonacci and Lucas Perfect Powers} 
\author{Yann Bugeaud, Maurice Mignotte, Samir Siksek}
\address{Yann Bugeaud\\
	Universit\'{e} Louis Pasteur\\
	U. F. R. de math\'{e}matiques\\
	7, rue Ren\'{e} Descartes\\
	67084 Strasbourg Cedex\\
	France}
\email{bugeaud@math.u-strasbg.fr}

\address{Maurice Mignotte\\
        Universit\'{e} Louis Pasteur\\
        U. F. R. de math\'{e}matiques\\
        7, rue Ren\'{e} Descartes\\
        67084 Strasbourg Cedex\\
        France}
\email{maurice@math.u-strasbg.fr}

\address{Samir Siksek\\
	Department of Mathematics and Statistics\\
        College  of Science \\
        Sultan Qaboos University \\
        P.O. Box 36 \\
        Al-Khod 123 \\
        Oman}
\email{siksek@squ.edu.om}

\date{\today}
\thanks{S. Siksek's work is funded by a grant from Sultan Qaboos University (IG/SCI/DOMS/02/06)}


\begin {abstract}
This is the first in a series of papers whereby we combine the
classical approach to exponential Diophantine equations (linear
forms in logarithms, Thue equations, etc.) with a modular 
approach based on some of the ideas of the proof of Fermat's Last 
Theorem. 
In this paper we give new improved bounds for linear forms
in three logarithms. We also apply a combination
of  classical techniques with the modular approach to
show that the only perfect powers in the Fibonacci
sequence are $0,1,8,144$ and the only perfect powers in the Lucas
sequence are $1,4$.  
\end {abstract}
\maketitle
\section{Introduction}

Wiles' proof of Fermat's Last Theorem \cite{Wi}, \cite{TW} is certainly
the most spectacular recent achievement in the field of Diophantine equations.
The proof uses what may be called the \lq modular\rq\ approach
which has since been applied to many other Diophantine equations; 
mostly -- though not exclusively -- of the form
\begin{equation}\label{eqn:Fermat}
a x^p+b y^p = c z^p, \qquad ax^p +b y^p = c z^2, \qquad ax^p+b y^p = c z^3, \ldots \quad \text{($p$ prime)}.
\end{equation}
The strategy  of the modular approach is simple enough: 
associate to a putative solution of such a Diophantine 
equation an elliptic curve, called a Frey curve, in a way that the discriminant is 
a $p$-th power up to some small factor. Next (under some technical assumptions) 
apply Ribet's level-lowering theorem \cite{Ribet}
to show that the Galois representation on the $p$-torsion of the Frey curve
arises from a newform
of weight $2$ and a fairly small level $N$ say. If there are no such newforms 
then there are no non-trivial solutions to the original Diophantine equation
(a solution is trivial if the corresponding Frey curve is singular).
Occasionally, even when one has newforms of the predicted level there
is still a possibility of somehow showing that it is incompatible with
the original Galois representation (see for example \cite{DM}, \cite {BS}), 
though there does not seem to be a general strategy.

A fact that has been underexploited is that the modular approach 
yields a tremendous amount of local information about the solutions
of the Diophantine equations. For equations of the form~(\ref{eqn:Fermat})
it is perhaps difficult to exploit this information successfully since we
neither know of a bound for the exponent $p$, nor for the variables $x,y,z$.
This suggests that the modular approach should be applied to exponential
Diophantine equations; for example, equations of the form
\[
a x^p + b y^p = c, \qquad a x^2+b=c y^p, \ldots \quad \text{($p$ prime)}.
\]
For such equations, Baker's theory of linear forms in logarithms (see the
book of Shorey and Tijdeman \cite{ST}) gives bounds for both the
exponent $p$ and the variables $x,y$.  This approach through linear forms
in logarithms and Thue equations, which we 
term the \lq classical\rq\ approach, has undergone through substantial
refinements, though it still often yields bounds that can only be described
as astronomical.

The present paper is the first in a series of papers whose aims are 
the following:
\begin{enumerate}
\item[(I)] To present theoretical improvements to various aspects of
the classical approach.
\item[(II)] To show how local information obtained through the modular
approach can be used to reduce the size of the bounds, both for 
exponents and for variables, of solutions to exponential Diophantine
equations.
\item[(III)] To show how local information obtained through the modular
approach can be pieced together to provide a proof that there are no
missing solutions less than the bounds obtained in (I),(II).
\item[(IV)] To solve various famous and hitherto outstanding exponential
Diophantine equations.
\end{enumerate}

Our theoretical improvement in this paper is a new and powerful
lower bound for linear forms in three logarithms. Such a lower
bound is often the key to bounding the exponent in an exponential
Diophantine equation. This is our
choice for (I). Our choice for (IV) is the infamous problem
of determining all perfect powers in the Fibonacci and Lucas
sequences. Items (II),(III) will be present in this paper
only in the context of solving this problem. A sequel 
combining the classical and modular approaches for 
Diophantine equations of the form $x^2 \pm D=y^p$
is in preparation~\cite{BMS}.

We delay presenting our lower bound for linear forms in three
logarithms till Section~\ref{sec:3logs}, as this is somewhat technical.
Regarding the Fibonacci and Lucas sequences we prove
the following Theorems.
\begin{thm}\label{thm:Fibonacci}
	Let $F_n$ be the $n$-th term of the Fibonacci
	sequence defined by $F_0=0, F_1=1$ and  $F_{n+2}=F_{n+1}+F_n$ for $n \geq 0$.
	The only perfect powers in the Fibonacci 
	sequence are $F_0=0$, $F_1=1$, $F_2=1$, $F_6=8$, $F_{12}=144$.
\end{thm}
\begin{thm}\label{thm:Lucas}
	Let $L_n$ be the $n$-th term of the Lucas 
        sequence defined by $L_0=2, L_1=1$ and  $L_{n+2}=L_{n+1}+L_n$ for $n \geq 0$.
	The only perfect powers in the Lucas
	sequence are $L_1=1$, $L_3=4$.
\end{thm}

It is appropriate to point out that equations
$F_n=y^p$ and $L_n=y^p$ have
previously been solved for small values of the exponent $p$
by various authors.
We present a brief survey of known results in Section~\ref{sec:survey}.

The main steps in the proofs of Theorems \ref{thm:Fibonacci} and \ref{thm:Lucas}
are as follows: 
\begin{enumerate}
	\item[(i)] We associate putative
	solutions to the equations 
	$F_n=y^p$ and $L_n=y^p$ with even index $n$ 
	to Frey curves and apply level-lowering.
	This, together with	
	some elementary arguments is used reduce to
	the case where the index $n$ satisfies 
	$n \equiv \pm 1 \pmod{6}$ 
	for equations $F_n=y^p$ and
	$L_n=y^p$.

	\item[(ii)] We then show that we may 
	suppose that the index $n$ 
	in the equations $F_n=y^p$ and $L_n=y^p$ 
	is prime. In the Fibonacci case this
	is essentially a result proved 
	first by Peth\H{o} \cite{Pe83} and Robbins \cite{Rob} (independently).
	 
	\item[(iii)] We apply level-lowering again
	under the assumption that the index
	$n$ is odd. We are able to show
	using this that $n \equiv \pm 1 \pmod{p}$
	for $p < 2 \times 10^{8}$ in the Fibonacci case.
	In the Lucas case we prove that
	$n \equiv \pm 1 \pmod{p}$ unconditionally.

	\item[(iv)] We show how to reduce 
	the equations $F_n=y^p$ and $L_n=y^p$
	to Thue equations.
	We do not solve these Thue equations completely, 
	but compute explicit upper bounds for their
	solutions using classical methods 
	(see for example \cite{BG}). 
	This provides us with upper bounds for $n$ in terms of $p$.
	In the Lucas case we need the fact that
	$n \equiv \pm 1 \pmod{p}$ to obtain 
	a simpler equation of Thue type.

	\item[(v)]   We show how the results of the level-lowering
	of step (iii) can be used,
	with the aid of a computer program, to
	produce extremely stringent congruence conditions on $n$. For 
	$p \leq 733$ in the Fibonacci case, and for $p \leq 281$
	in the Lucas case, the congruences obtained are so strong
	that, when combined with the upper bounds for $n$ in terms
	of $p$ obtained in (iv), give a complete resolution for $F_n=y^p$
	and $L_n=y^p$.

	\item[(vi)] It is known that the equation $L_n=y^p$
	yields a linear form in two logarithms. Applying 
	the bounds of Laurent, Mignotte and Nesterenko \cite{LMN} 
	we show that $p \leq 281$ in the Lucas case. This completes
	the determination of perfect powers in the Lucas sequences.

	\item[(vii)] The equation $F_n=y^p$ yields a linear form
	in three logarithms. However if $p< 2 \times {10}^8$ then by step (iii)
	we know that $n \equiv \pm 1 \pmod{p}$. We show how
	in this case the linear form in three logarithms may
	be rewritten as a linear form in two logarithms.
	Applying \cite{LMN} we deduce that $p \leq 733$
	which we have already solved in step (v). 

	\item[(viii)] To complete the resolution of $F_n=y^p$
	it is enough to show that $p < 2 \times {10}^8$. We present
	a powerful improvement to known bounds for linear forms
	in $3$ logarithms. Applying our result shows indeed
	that $p < 2 \times {10}^8$ and this completes the determination
	of perfect powers in the Fibonacci sequence.
\end{enumerate}	

Let us make some brief comments.

The condition $n \equiv \pm 1 \pmod{p}$ obtained after step (iii)
cannot be strengthened. Indeed, we may define $F_n$ and $L_n$ for negative
$n$ by the recursion formulae $F_{n+2} = F_{n+1} + F_n$ and
$L_{n+2} = L_{n+1} + L_n$. We then observe that $F_{-1} = 1$
and $L_{-1} = -1$. Consequently, $F_{-1}$, $F_1$, $L_{-1}$ and
$L_1$ are $p$-th powers for any odd prime $p$. Thus equations
$F_n = y^p$ and $L_n = y^p$ do have solutions with
$n \equiv \pm 1 \pmod{p}$.

The strategy of combining explicit upper bounds for the
solutions of Thue equations with a sieve has already been applied
successfully in \cite{BMRS}. The idea of combining explicit
upper bounds with the modular approach was first tentatively
floated in \cite{SC}. 

A crucial observation for the proof of Theorem 1 is the fact that, with a
modicum of computation, we can indeed use linear forms in two
logarithms, and then get a much smaller upper bound for the
exponent $p$.

The present paper is organised as follows. Section 2 is devoted to a 
survey of previous results. Sections 3 and 4 are concerned with
useful preliminaries. Steps (i) and (ii) are treated in Sections 5 and 6,
respectively. Sections 7 and 8 are devoted to step (iii). Sections 9
and 10 are concerned with Steps (iv) and (v). Section 11 deals with 
steps (vi) and (vii), and finishes the proof of Theorem 2. Finally,
the proof of Theorem 1 is completed in Section 13, which deals with
step (viii), by applying estimates for linear forms in three
logarithms proved in Section 12.

The computations in the paper were performed using the computer
packages {\tt PARI/GP} \cite{PARI} and {\tt MAGMA} \cite{Magma}.
The total running time for the various computational parts of 
the proof of Theorem~\ref{thm:Fibonacci} is roughly 158 hours on
a $1.7$ GHz Intel Pentium 4. By contrast, the total time for 
the corresponding computational parts of the proof of 
Theorem~\ref{thm:Lucas} is roughly 6 hours.

\section{A Brief Survey of Previous Results}\label{sec:survey}
In this section we would like place our 
Theorems~\ref{thm:Fibonacci} and~\ref{thm:Lucas} 
in the context of other exponential
Diophantine equations. We also give a very brief survey
of results known to us on the problem of perfect powers in the Fibonacci
and Lucas sequences, though we make no claim that our survey is exhaustive.

Thanks to Baker's theory of linear forms in logarithms, we know
(see for example the book of Shorey and Tijdeman \cite{ST}) that many families of
Diophantine equations have finitely many integer solutions, and that one can even
compute upper bounds for their absolute values. These upper bounds
are however huge and do not enable us to provide complete
lists of solutions by brutal enumeration. During the last decade,
thanks to important progress in computational number
theory (such as the LLL-algorithm) and also in the theory of linear forms
in logarithms (the numerical constants have been substantially reduced
in comparison to Baker's first papers), we are now able to solve
completely some exponential Diophantine equations. Perhaps the most
striking achievement obtained via techniques from
Diophantine approximation is a result of Bennett \cite{Be}, asserting that,
for any integers $a$, $b$ and $p \geq 3$ with $a > b \geq 1$,
the Diophantine equation
$$
|a X^p - b Y^p| = 1
$$
has at most one solution in positive integers $X$ and $Y$.

Among other results in this area
obtained thanks to (at least in part) to theory of linear forms
in logarithms, let us
quote that Bugeaud and Mignotte \cite{BM} proved that the equation
$(10^n - 1) / (10- 1) = y^p$ has no solution with $y>1$, and that
Bilu, Hanrot and Voutier \cite{BHV} solved the long-standing problem
of the existence of primitive divisors of Lucas--Lehmer sequences.

Despite substantial theoretical progress and the use of techniques coming
from arithmetic geometry and developed in connection with Fermat's
Last Theorem (see for example the paper of Bennett and Skinner \cite{BS}),
some celebrated Diophantine equations are still unsolved. We would particularly like
to draw the reader's attention to the following three equations:
\begin{equation}\label{eqn:SC}
x^2 + 7 = y^p,  \qquad p \geq 3,
\end{equation}
\begin{equation}\label{eqn:BS}
x^2 - 2 = y^p,  \qquad p \geq 3,
\end{equation}
and
\begin{equation}\label{eqn:Fibytm}
F_n = y^p,  \qquad \text{$n \geq 0$ and $p \geq 2$},
\end{equation}
where $F_n$ is the $n$-th term in the Fibonacci sequence.
Let us explain the difficulties encountered with equations (\ref{eqn:SC}),
(\ref{eqn:BS}) and (\ref{eqn:Fibytm}). Classically, we first use estimates
for linear forms in logarithms
in order to bound the exponent $p$, and then we perform a sieve.
Equations (\ref{eqn:SC}) and (\ref{eqn:Fibytm}) yield linear forms in
three logarithms, and thus
upper bounds for $p$ of the order of $10^{13}$;
at present far too large to allow the complete resolution of (\ref{eqn:SC})
and (\ref{eqn:Fibytm}) using classical methods (however, a promising attempt at
equation~(\ref{eqn:SC}) is made in \cite{SC}). The case of (\ref{eqn:BS}) is
different, since estimates for linear forms in two logarithms
yield that $n$ is at most $164969$ \cite{Iv}, an upper bound which
can certainly be (at least) slightly improved.
There is however a notorious difficulty in (\ref{eqn:BS}) and (\ref{eqn:Fibytm}),
namely the existence
of solutions $1^2 - 2 = (-1)^p$ and $F_1 = 1^p$ for each value of the exponent
$p$.  These small solutions prevent us from using a sieve as efficient as
the one used for (\ref{eqn:SC}).
A natural way to overcome this
is to derive from (\ref{eqn:BS}) and (\ref{eqn:Fibytm}) Thue equations,
though these are of degree far too large to allow for a complete resolution
using classical methods alone.

As we have explained in the introduction, the present work is devoted to 
equation~(\ref{eqn:Fibytm}), and to the analogous equation for the Lucas sequence. 

As for general results, Peth\H o \cite{Pe82} and, independently,
Shorey and Stewart \cite{ShSt} proved that there are only finitely many
perfect powers in any non-trivial binary recurrence sequence. Their
proofs, based on Baker's theory of linear forms in logarithms, 
are effective but yield huge bounds. We now turn to specific results
on the Fibonacci and Lucas sequences.

\begin{itemize}
\item The only perfect squares in the Fibonacci sequence are 
$F_0=0,~F_1=F_2=1,~F_{12}=144$;
 see \cite{Cohn1} and \cite{Wy}.
\item The only perfect cubes in the Fibonacci sequence are 
$F_0=0,~F_1=F_2=1,~F_6=8$;
see \cite{LF}.
\item  For $m=5,7,11,13,17$, the only $m$-th powers are $F_0=0,~F_1=F_2=1$.
This was proved by J. McLaughlin \cite{Mc} by using a linear form in logarithms
together with the LLL algorithm.
\item If $n > 2 $ and $F_n=y^p$  
then $p < 5.1 \times 10^{17}$; this was proved by
Peth\H{o} using a linear form in three logarithms \cite{Pe01}.
In the same paper he also showed that if $n>2$ and $L_n=y^p$ 
then $p<13222$ using a linear form in two logarithms.
\item Another result which is particularly relevant to us is the 
following: If $p \geq 3$ and $F_n=y^p$ for integer $y$ then either $n=0,1,2,6$ or
there is a prime $q \mid n$ such that $F_q=y_1^p$, for some integer $y_1$.
This result was established by Peth\H{o} \cite{Pe83} and Robbins \cite{Rob} independently.
\item Cohn \cite{Cohn3} proved that $L_1 = 1$ and
$L_3 = 4$ are the only squares in the Lucas sequence. 
\item London and Finkelstein \cite{LF} proved that $L_1 = 1$ is the 
only cube in the Lucas sequence.
\end{itemize}

\section{Preliminaries}
We collect in this section various results
which will be useful throughout this paper.
Our problem of determining the
perfect powers in the Fibonacci and Lucas sequences
naturally reduces to the problem of solving the
following pair of equations:
\begin{equation}\label{eqn:Fibytp}
	F_n=y^p, \qquad \text{$n \geq 0$, and $p$ prime,}
\end{equation}
and
\begin{equation}\label{eqn:Lucytp}
	L_n=y^p, \qquad \text{$n \geq 0$, and $p$ prime}.
\end{equation}
Throughout this paper we will use
the fact that
\begin{equation}\label{eqn:FL}
	F_n =\frac{\omega^n - \tau^n}{\sqrt{5}},
	\qquad L_n = \omega^n + \tau^n,
\end{equation}
where
\begin{equation}\label{eqn:lmu}
	\omega = \frac{1+\sqrt{5}}{2},
	\qquad \tau = \frac{1- \sqrt{5}}{2}.
\end{equation}
This quickly leads us to associate
the equations $F_n=y^p$ and $L_n=y^p$
with auxiliary equations as the following
two Lemmas show.
\begin{lem}\label{lem:FibAux}
	Suppose that $F_n=y^p$. If $n$ is odd then 
		\begin{equation}\label{eqn:Fibodd}
			5y^{2p}=L_n^2+4,	
		\end{equation}
	and if $n$ is even then
		\begin{equation}\label{eqn:Fibeven}
			5y^{2p}=L_n^2-4.	
		\end{equation}
\end{lem}
\begin{lem}\label{lem:LucAux}
        Suppose that $L_n=y^p$. If $n$ is odd then
                \begin{equation}\label{eqn:Lucodd}
        		y^{2p}=5 F_n^2-4,
	        \end{equation}
        and if $n$ is even then
                \begin{equation}\label{eqn:Luceven}
        		y^{2p}=5 F_n^2+4.
	        \end{equation}
\end{lem}

For a prime $l \neq 5$ define
\begin{equation}\label{eqn:M(l)}
	M(l)=
		\left\{
		\begin{array}{ll}
			l-1 & \text{if $l \equiv \pm 1 \pmod{5}$} \\
			2(l+1) & \text{if $l \equiv \pm 2 \pmod{5}$.}
		\end{array}
		\right.
\end{equation}
We will need the following two Lemmas.
\begin{lem}\label{lem:M(l)}
Suppose that $l \neq 5$ is a prime and $n \equiv m \pmod{M(l)}$.
Then $F_n \equiv F_m \pmod{l}$ and $L_n \equiv L_m \pmod{l}$.
\end{lem}
\begin{proof}
	Write $\OO$ for the ring of integers of the field
	$\Q(\sqrt{5})$.
	Recall, by~(\ref{eqn:FL}), that $F_n$ and $L_n$ are 
	expressed in terms of
	$\omega,\tau$. Let $\pi$ be a prime in $\OO$
	dividing $l$. To prove the Lemma all we need to show
	is that
	\[
		\omega^{M(l)} \equiv \tau^{M(l)} \equiv 1 \pmod{\pi}.
	\]
	If $l \equiv \pm 1 \pmod{5}$ then
	$5$ is a quadratic residue modulo $l$.
	The Lemma follows immediately in this 
	case from the fact that 
	$(\OO/{\pi \OO})^* \cong \F_l^*$ and
	so has order $l-1$. 

	Now suppose that $l \equiv \pm 2 \pmod{5}$.
	Note that
	\[
		\omega^l \equiv \frac{1^l + 5^\frac{l-1}{2} \sqrt{5}}{2^l}
		\equiv \frac{1-\sqrt{5}}{2} \equiv \tau \pmod{\pi},
	\]
	since $5$ is a quadratic non-residue modulo $l$. Thus
	\[
		\omega^{M(l)} \equiv \omega^{2(l+1)} \equiv 
		(\omega \tau)^2 \equiv 1 \pmod{\pi},
	\]
	and similarly for $\tau$.	
\end{proof}

\begin{lem}\label{lem:cong}
	The residues of $L_n,~F_n$ modulo $4$ depend only on the
	residue of $n$ modulo $6$, and are given by the following table
	\begin{center}
	\begin{tabular}{c || c | c }
        		&       $L_n \pmod{4}$  &       $F_n \pmod{4}$  \\
		\hline
		\hline
		$n \equiv 0 \pmod{6}$   &       $2$             &               $0$     \\
		$n \equiv 1 \pmod{6}$     &       $1$               &               $1$ \\
		$n \equiv 2 \pmod{6}$     &       $3$               &               $1$ \\
		$n \equiv 3 \pmod{6}$     &       $0$               &               $2$ \\
		$n \equiv 4 \pmod{6}$     &       $3$               &               $3$ \\
		$n \equiv 5 \pmod{6}$     &       $3$               &               $1$
	\end{tabular}
	\end{center}
\end{lem}
\begin{proof}
	The Lemma is proved by a straightforward induction, using the 
recurrence relations defining $F_n$ and $L_n$.
\end{proof}


\section{Eliminating Small Exponents and Indices}
We will later need to assume that the exponent $p$
and the index $n$ in the equations~(\ref{eqn:Fibytp})
and~(\ref{eqn:Lucytp}) are not too small. More precisely,
in this Section, we prove the following pair of Propositions.
\begin{prop}\label{prop:FibonacciSimplified}
	If there is a perfect power in the Fibonacci sequence
	not listed in Theorem~\ref{thm:Fibonacci}
	then there is a solution to the equation
	\begin{equation}\label{eqn:FibonacciSimplified}
	        F_n=y^p, \qquad \text{$n> 25000$ and $p \geq 7$ is prime.}
	\end{equation}
\end{prop}
\begin{prop}\label{prop:LucasSimplified}
        If there is a perfect power in the Lucas sequence
        not listed in Theorem~\ref{thm:Lucas}
        then there is a solution to the equation
	\begin{equation}\label{eqn:LucasSimplified}
	        L_n=y^p, \qquad \text{$n> 25000$  and $p \geq 7$ is prime.}
	\end{equation}
\end{prop}
The propositions follow from the results on Fibonacci perfect powers
quoted in Section~\ref{sec:survey} together with Lemmas~\ref{lem:nlarge}
and~\ref{lem:Luc235} below.
\subsection{Ruling Out Small Values of the Index $n$}
\begin{lem}\label{lem:nlarge}
	For no integer $13 \leq n \leq 25000$ is $F_n$ a perfect power.
	For no integer $4 \leq n \leq 25000$ is $L_n$ a perfect power.
\end{lem}
\begin{proof}
	Suppose $F_n=y^p$ where $p$ is some prime and $n$
	is in the range $13 \leq n \leq 25000$. It is easy to see
	from~(\ref{eqn:FL}), (\ref{eqn:lmu}) that 
	$2 \leq p \leq n \log(\omega)/\log(2)$. Now fix $n,~p$
	and we would like to show that $F_n$ is not a $p$-th power.

	Suppose $l$ is a prime satisfying $l \equiv \pm 1 \pmod{5}$ and
	$l \equiv 1 \pmod{p}$. The condition $l \equiv \pm 1 \pmod{5}$
	ensures that $5$ is a quadratic residue modulo $l$.
	Then one can easily compute $F_n$ modulo $l$ using~(\ref{eqn:FL}) 
	(without having to write down
	$F_n$). Now let $k=(l-1)/p$. If $F_n^k \not \equiv 1 \pmod{l}$
	then we know that $F_n$ is not a $p$-th power.

	We wrote a short {\tt PARI/GP} program to check for $n$ in the above
	range, and for each prime $2 \leq p \leq n \log(\omega)/\log(2)$
	that there exits a prime $l$ proving that $F_n$ is not a $p$-th power,
	using the above idea.
	This took roughly 15 minutes on a 1.7 GHz Pentium 4.

	The corresponding result for the Lucas sequence is proved 
	in exactly the same way, with program taking roughly 16 minutes to run
	on the same machine.
\end{proof}

\subsection{Solutions with Exponent $p=2,3,5$}
Later on when we come to apply level-lowering
we will need to assume that $p \geq 7$.
It is  straightforward to solve equations~(\ref{eqn:Fibytp})
and~(\ref{eqn:Lucytp}) for $p=2,~3,~5$
with the help of the computer algebra package
{\tt Magma}. We give the details for the Lucas case;
the Fibonacci case is similar. Alternatively we
could quote the known results surveyed in Section~\ref{sec:survey},
although $p=5$ for the Lucas case does not seem to be covered by the
literature.

\begin{lem}\label{lem:Luc235}
	The only solutions to the equation~(\ref{eqn:Lucytp}) with $p=2,~3,~5$
	are $(n,y,p)=(1,1,p)$ and $(3,2,2)$.
\end{lem}
\begin{proof}
	Suppose first that $n$ is even. By Lemma~\ref{lem:LucAux}
	it is enough to show that~(\ref{eqn:Luceven}) does not have a solution.
	Suppose that $(n,y,p)$ is a solution to~(\ref{eqn:Luceven}). Clearly
	$F_n$ and $y$ are odd, and $y$ is not divisible by $5$. 
	Thus we have
	\[
		(2+F_n \sqrt{-5})= \mathfrak{a}^{2p}
	\]
	for some ideal $\mathfrak{a}$ of $\Z[\sqrt{-5}]$. 
	Now the class number of $\Z[\sqrt{-5}]$
	is $2$, and hence $\mathfrak{a}^2$ is a principal ideal. 
	It follows that
	\[
		2+F_n \sqrt{-5}=(u+v\sqrt{-5})^p
	\]
	for some integers $u,~v$. If $p=2$ then we get $2=u^2-5v^2$ which
	is impossible modulo $5$. If $p=3$ then
	\[
		2=u(u^2-15v^2),
	\]
	and if $p=5$ then
	\[
		2=u(u^4-50u^2 v^2+125 v^4).
	\]
	It is straightforward to see that both of these are impossible.
	Next we turn to the case where $n$ is odd. Again by Lemma~\ref{lem:LucAux} 
	it is enough to solve the equation~(\ref{eqn:Lucodd}). Suppose first that
	$p=3,5$. If $(n,y,p)$ is any solution to equation~(\ref{eqn:Lucodd}) then we 
	quickly see that $y$ must be odd and
	\[
		2+\sqrt{5}F_n=\left( \frac{1+\sqrt{5}}{2} \right)^r (u+v \sqrt{5})^p
	\]
	For some $r=0,\ldots, p-1$ and $u,v$ are both integers or both halves of odd
	integers. The computer algebra package {\tt Magma} quickly solves
	all the resulting Thue equations showing that $y=\pm 1$. This implies that for
	$p=3,5$ the only solution to equation~(\ref{eqn:Lucytp}) is $(1,1,p)$.
	
	Finally to deal with $p=2$ we note that if $(n,y)$ satisfies~(\ref{eqn:Lucodd})
	then $(X,Y)=(5y^2,25 F_n y)$ is an integral point on the elliptic curve 
	$Y^2=X^3+100X$. Again {\tt Magma} quickly computes all integral points on this
	curve:
	these are $(X,Y)=(0,0),(5, \pm 25), (20,\pm 100)$, which yields the solutions
	$(n,y)=(1,1),(3,2)$. This completes the proof of the Lemma.
\end{proof}

\section{Reducing to the Case $n \equiv \pm 1 \pmod{6}$}
In this section we would like to reduce
the study of equations~(\ref{eqn:Fibytp})
and~(\ref{eqn:Lucytp}) to the special
case where the index $n$ satisfies $n \equiv \pm 1 \pmod{6}$. 
For Fibonacci we show that if there is some solution 
$(n,y,p)$ to~(\ref{eqn:Fibytp}) then there is 
another solution with the same exponent $p$ such that
the index $n$ satisfies the above condition.  
For Lucas sequences we prove 
the following stronger result. 
\begin{lem}\label{lem:LucCong}
        If $(n,y,p)$ is a solution to the equation~(\ref{eqn:Lucytp}) with $p \geq 7$
        then $n \equiv \pm 1 \pmod{6}$.
\end{lem}
For Fibonacci our result is weaker but still
useful.
\begin{lem}\label{lem:FibCong}
	If $(n,y,p)$ is a solution to equation~(\ref{eqn:Fibytp}) with $p \geq 7$
	then either $n=0$ or $n \equiv \pm 1 \pmod{6}$ or else $n=2k$ with 
	\begin{enumerate}
		\item[(a)] $k \equiv \pm 1 \pmod{6}$
		\item[(b)] $F_{k}=U^p$ and $L_{k}=V^p$ for some positive integers
		$U$ and $V$.
	\end{enumerate}
\end{lem}
The proofs of both of Lemmas~\ref{lem:LucCong} and~\ref{lem:FibCong}
make use of Frey curves and level-lowering. Here and elsewhere
where we make use these tools, we do not directly apply the original
results in this field (Ribet's level-lowering Theorem \cite{Ribet}, modularity
of elliptic curves by Wiles and others \cite{Wi}, \cite{Mod}, irreducibility of Galois 
representations by Mazur and others \cite{Mazur}, etc.). We will instead quote
directly from excellent recent paper of Bennett and Skinner \cite{BS},
which is concerned with equations of the form $Ax^n+By^n=Cz^2$. 
In every instance we will put our equation in this form before applying
the results of \cite{BS}.

\begin{proof}[Proof of Lemma~\ref{lem:LucCong}]
	Suppose that $(n,y,p)$ is a solution to equation~(\ref{eqn:Lucytp})
	with  $p \geq 7$.

	We observe first that $n \not \equiv 0,3 \pmod{6}$.
	For in this case Lemma~\ref{lem:cong} implies that both
	$F_n$ and $L_n$ are even, and hence by Lemma~\ref{lem:LucAux} 
	either $5$ or $-5$ is a $2$-adic square, which is not the case.

	We now restrict our attention to $n \equiv 2,~4 \pmod{6}$ and $p \geq 7$
	and show that this leads to a contradiction. This is enough to 
	prove the Lemma. Let
	\begin{eqnarray*}
        	G_n=\left\{
        	\begin{array}{ll}
                	-F_n & \text{if\ } n \equiv 2 \pmod{6} \\
                	F_n & \text{if\ } n \equiv 4 \pmod{6}.
        	\end{array}
        	\right.
	\end{eqnarray*}
	It follows from Lemma~\ref{lem:LucAux} that
	\[
        	y^{2p}=5 {G_n}^2 +4.
	\]
	We associate to our solution $(n,y,p)$ of~(\ref{eqn:Lucytp}) with 
	$n \equiv 2,~4 \pmod{6}$ the Frey curve
	\begin{equation}
        	E_{n}: \qquad Y^2 = X^3+ 5 G_n X^2- 5 X.
	\end{equation}
	Let $E$ be the elliptic curve 100A1
	in Cremona's tables \cite{Cre}; $E$ has the following model:
	\[
        	E: \quad Y^2=X^3-X^2-33X+62.
	\]
	Write $\rho_p(E)$ for the Galois representation
	\[
        	\rho_p(E) : {\rm Gal}(\overline{\Q}/\Q) \rightarrow {\rm Aut}(E[p])
	\]
	on the $p$-torsion of $E$, and let $\rho_p(E_n)$ be the corresponding
	Galois representation for $E_n$.

	Applying the results of  \cite[Sections 2,~3]{BS},
	we see that $\rho_p(E_n)$ arises from a cuspidal newform of
	weight $2$,  level $100$, and trivial Nebentypus character. However using
	the computer algebra package {\tt Magma}  
	we find that the dimension of newforms of
	weight $2$ and level $100$
	is $1$. Moreover the curve $E$ above is (up to isogeny) the unique elliptic curve
	of conductor $100$. Thus $\rho_p(E_n)$ and $\rho_p(E)$ are isomorphic. It follows
	from this, by \cite[Proposition 4.4]{BS}, that $5$ does not divide the
	denominator of the $j$-invariant of $E$. This is not true as
	$j(E)=16384/5$, giving us a contradiction.

	For the convenience of the reader we point out that in Bennett and Skinner's
	notation:
	\[
	A=1, \quad B=-4, \quad C=5, \quad a=y^2, \quad b=1, \quad c=G_n.
	\]
	Lemma~\ref{lem:cong} and our definition of $G_n$ above
	imply that $c \equiv 3 \pmod{4}$ which is needed
	to apply the results of Bennett and Skinner.
	This completes our proof of Lemma~\ref{lem:LucCong}.
\end{proof}

\begin{proof}[Proof of Lemma~\ref{lem:FibCong}]
	Suppose that $(n,y,p)$ is a solution to 
	equation~(\ref{eqn:Fibytp}) with $n \neq 0$ and $p \geq 7$. 
	By Lemma~\ref{lem:cong}
	we see that $n \not \equiv 3 \pmod{6}$. 
	Suppose then that $n \not \equiv \pm 1 \pmod{6}$.
	Clearly $n=2k$ for some integer $k$.
	It is well-known and easy to see from~(\ref{eqn:FL})  
	that $F_{n}=F_{2k}=F_k L_k$.
	It is also easy to see that the greatest common divisor of 
	$F_k$ and $L_k$ is either $1$ or $2$. 
	The crux of the proof is to show that if $F_n=y^p$ then $F_n$ is odd.

	Thus suppose that $F_n$ (and hence $y$) is even.
	Lemma~\ref{lem:LucAux} tells us that $5 y^{2p}+4=L_{n}^2$.
	Since $y$ even, we see that  $2 \mid \mid L_{n}$.
	Let $z=y/2$ and
	\[
		x = \left\{
		\begin{array}{ll}
			L_{n}/2 & \text{if $L_{n} \equiv 2 \pmod{8}$} \\
			-L_{n}/2 & \text{if $L_{n} \equiv 6 \pmod{8}$ .}
		\end{array}
		\right.
	\]

	Thus $x \equiv 1 \pmod{4}$ and 
	\[
	2^{2p-2}\cdot 5 z^{2p}+1=x^2.
	\]
	Following \cite[Section 2]{BS} we associate to this equation
	the Frey curve
	\[
	Y^2+XY=X^3+\left(\frac{x-1}{4}\right)X^2+ 2^{2p-8} \cdot 5 z^{2p} X.
	\]
	 Applying level-lowering \cite[Section 3]{BS} shows that the Galois
	representation arises from a cusp form of weight $2$ and level $10$.
	Since there are no such cusp forms we get a contradiction (This is 
	essentially the same argument used in the proof of Fermat's Last
	Theorem). It is noted that the argument here fails for $n=0$
	since in this case the Frey curve is singular.

	We deduce that $F_{n}$ is odd, and so that $F_k=U^p$ and $L_k=V^p$
	for some positive integers $U,~V$. 
	By Lemma~\ref{lem:LucCong} we know that $k \equiv \pm 1 \pmod{6}$.
	This completes the proof of Lemma~\ref{lem:FibCong}.
\end{proof}
\section{Reduction to the Prime Index Case}
In this section we reduce our problem to the assumption that
the index $n$ is prime, as in the following pair of Propositions. 
\begin{prop}\label{prop:FibPrimeIndex}
        If there is a perfect power in the Fibonacci sequence
        not listed in Theorem~\ref{thm:Fibonacci}
        then there is a solution to the equation
        \begin{equation}\label{eqn:FibPrimeIndex}
                F_n=y^p, \qquad \text{$n> 25000,~p \geq 7$ with $n,~p$ prime.}
        \end{equation}
\end{prop}
\begin{prop}\label{prop:LucPrimeIndex}
        If there is a perfect power in the Lucas sequence
        not listed in Theorem~\ref{thm:Lucas}
        then there is a solution to the equation
        \begin{equation}\label{eqn:LucPrimeIndex}
                L_n=y^p, \qquad \text{$n> 25000,~p \geq 7$ with $n,~p$ prime.}
        \end{equation}
\end{prop}
After proving these two propositions the remainder of this paper
will be devoted to showing the there are no solutions to equations
~(\ref{eqn:FibPrimeIndex}) and~(\ref{eqn:LucPrimeIndex}).
\begin{proof}[Proof of Proposition~\ref{prop:FibPrimeIndex}]
	If $F_n=y^p$ with $n$ is odd then this just the result of
	Peth\H{o} and Robbins quoted
	in Section~\ref{sec:survey} together with our 
	Proposition~\ref{prop:FibonacciSimplified}.
	Suppose $n=2k$. By Lemma~\ref{lem:FibCong} 
	we know that $k$ is odd and $F_k=U^p$ for some integer
	$U$. Now simply apply the result of Peth\H{o} 
	and Robbins again together
	with Proposition~\ref{prop:FibonacciSimplified}.
\end{proof}
\begin{proof}[Proof of Proposition~\ref{prop:LucPrimeIndex}]
	Suppose that $L_n=y^p$ where $n \neq 1,~3$ and $p \geq 7$.
	By Lemma~\ref{lem:LucCong} we know that $n \equiv \pm 1 \pmod{6}$,
	and so $n$ is odd. If $n$ is prime then the result follows
	from Proposition~\ref{prop:LucasSimplified}. Thus
	suppose that $n$ is
	composite and let $q$ be its smallest prime factor. Write $n=k q$,
	where $k > 1$.
	Then $L_n=y^p$ can be rewritten as
	\begin{equation}\label{eqn:factor}
		(\omega^{k} - \omega^{-k})(\omega^{k(q-1)}+\cdots+\omega^{-k(q-1)})
		=\omega^n-\omega^{-n} = y^p.
	\end{equation}
	It is straightforward to see that the two factors on the left-hand side are
	in $\Z$ and that their greatest common factor divides $q$.

	Suppose that $q$ divides the two factors. Then we see that
	\[
		\omega^{2k} \equiv 1 \pmod {\pi}
	\]
	for some prime $\pi$ of $\OO$ lying above $q$.
	But $\omega^{2}-1=\omega $ and so
	$\omega^{2} \not \equiv 1 \pmod{\pi}$.
	It follows that the order of the image of $\omega^2$
	in $(\OO/\pi)^*$ is not $1$ and that it divides $k$ and hence $n$.
	But $\# (\OO/\pi)^*$ is either $q-1$ or $q^2-1$.
	It follows that some non-trivial factor of $n$ divides
	$(q-1)(q+1)$. 
	Moreover, $n$ is odd, and all odd prime factors 
	of $(q-1)(q+1)$ are smaller than $q$.
	This contradicts the assumption that
	$q$ is the smallest prime factor of $n$.

	We deduce that $q$ does not divide the factors on the
	right-hand side of~(\ref{eqn:factor}). Hence
	$L_{k}=\omega^{k} - \omega^{-k}=y_1^p$ for some
	integer $y_1$. If $k$ is prime then the proof is complete
	by invoking Proposition~\ref{prop:LucasSimplified}. 
	Otherwise apply the above argument recursively.
\end{proof}

\section{Level-Lowering for Fibonacci - The Odd Index Case}\label{sec:FibLL}

Previously we have used a Frey curve and level-lowering to
obtain information about solutions of $F_n=y^p$ for even $n$.
In this section we associate a Frey curve to any
solution of equation~(\ref{eqn:FibPrimeIndex}).

Suppose that $(n,y,p)$ is a solution to~(\ref{eqn:FibPrimeIndex}).
Thus $n,~p$ are primes with $p \geq 7$ and $n > 25000$.
Let
\begin{equation}\label{eqn:Hn}
	H_n= \left\{
		\begin{array}{ll}
			L_n  & \text{ if } n \equiv 1 \pmod{6} \\
			-L_n  & \text{ if } n \equiv 5 \pmod{6}.
		\end{array}
	\right.
\end{equation}
\begin{lem}\label{lem:1mod4}
	With notation as above, $H_n \equiv 1 \pmod{4}$
	and
	\begin{equation}\label{eqn:hn2}
		5 y^{2p} -4 =H_n^2.
	\end{equation}
\end{lem}
The Lemma follows immediately from Lemma~\ref{lem:FibAux} and
Lemma~\ref{lem:cong}.

We associate to the solution $(n,y,p)$ the Frey curve
\begin{equation} \label{eqn:FibFrey}
	E_{n}: \qquad Y^2 = X^3+ H_n X^2 - X.
\end{equation}

We now come to level-lowering. Let $E$ be the following elliptic curve
over~$\Q$:
\begin{equation}\label{eqn:FibEll}
	E: \quad Y^2=X^3+X^2-X;
\end{equation}
this is curve 20A2 in Cremona's tables \cite{Cre}. As before, 
write $\rho_p(E)$ for the Galois representation
on the $p$-torsion of $E$, and let $\rho_p(E_n)$ be the corresponding Galois
representation on the $p$-torsion of $E_n$. If $l$ is a prime, let
$a_l(E)$ be the trace of Frobenius of the curve $E$ at $l$,
and let $a_l(E_n)$ denote the corresponding trace of Frobenius of $E_n$.

\begin{prop}\label{prop:FibLL}
Suppose that $(n,y,p)$ is a solution to~(\ref{eqn:FibPrimeIndex}).
With notation as above, the Galois representations $\rho_p(E_n),~\rho_p(E)$
are isomorphic. Moreover, for any prime $l \neq 2,~5$,
\begin{itemize}
	\item[(i)] $a_l(E_n) \equiv a_l(E) \pmod{p}$ if $ l \nmid y$. 
	\item[(ii)] $l+1 \equiv \pm a_l(E) \pmod{p}$
	 if $l \mid y$.
\end{itemize}
\end{prop}
\begin{proof}
First we apply the results of \cite[Sections 2,3]{BS}.  From 
these we know that $\rho_p(E_n)$ arises from a cuspidal newform
of weight $2$, level $20$, and trivial Nebentypus character.
(In applying the results of \cite{BS} we need Lemma~\ref{lem:1mod4}).
However $S_2(\Gamma_0(20))$ has dimension
$1$. Moreover, the curve $E$ is (up to isogeny) the unique curve of
conductor $20$. It follows that $\rho_p(E)$ and $\rho_p(E_n)$ are
isomorphic.

The rest of the Proposition follows from \cite[Proposition 3]{KO}, and the fact that if
$l \neq 2,5$ and $l \mid y$ then $l$ is a prime of multiplicative reduction
for $E_n$ and so $a_l(E_n)= \pm 1$.
\end{proof}

Proposition~\ref{prop:FibLL} is useful in several stages of our proof of
Theorem~\ref{thm:Fibonacci}. The following Proposition is needed
later, and follows from Proposition~\ref{prop:FibLL} and some computational
work.
\begin{prop}\label{prop:Fibnequivpm1}
	If $(n,y,p)$ is any solution to equation~(\ref{eqn:FibPrimeIndex})
	with $p < 2 \times {10}^8$ then $n \equiv \pm 1 \pmod{p}$.
\end{prop}

The idea behind the proof is inspired by a method of
Kraus (see \cite{Kra} or \cite{SC}) but there are added
complications in our situation: for any prime $p$
the equation $F_n=y^p$ has the solution $(n,y)=(1,1)$,
and also the solution $(n,y)=(-1,1)$ (got by extrapolating
the definition of the Fibonacci sequence backwards).

Before proving Proposition~\ref{prop:Fibnequivpm1} we start with
a little motivation. Suppose that $p \geq 7$
is a prime, and we find some small positive integer $k$
such that $l=2kp+1$ is prime, and $l \equiv \pm 1 \pmod{5}$.
It follows that $5$ is a quadratic residue modulo $l$,
and we choose an element in $\F_l$ which we conveniently
denote by $\sqrt{5}$, satisfying $(\sqrt{5})^2 \equiv 5 \pmod{l}$.
We may then think of $\omega,\tau$ (defined in~(\ref{eqn:lmu})) as elements
of $\F_l$.

Consider the equation $F_n=y^p$. Now $l-1=2kp$, with $k$ is
small. This means that $y^p$ comes from a small subset of $\F_l$.
We can now use the level-lowering to
predict the values of $y^p$. Hopefully, we may find that the
only value of $y^p$ modulo $l$ predicted by the level-lowering
and also belonging to our small subset are $\pm 1$. Under a further
minor hypothesis we can show that this implies that $n \equiv
\pm 1 \pmod{p}$. If a particular value of $k$ does not work we
may continue trying until a suitable $k$ is found.

We make all this precise.  Suppose as above that $l,~p$
are primes with $l=2kp+1$ and $l \equiv \pm 1 \pmod{5}$. Define
\[
	A(p,k)= \left\{
        	\zeta \in \left( \F_l^* \right)^{2p} \backslash \{ 1 \}
		\, : \quad \kro{5 \zeta -4}{l}= \text{$0$ or $1$}
        	\right\}.
\]
For each $\zeta \in A(p,k)$, choose an integer $\delta_\zeta$ such
that
\[
	\delta_\zeta^2 \equiv 5 \zeta -4 \pmod{l}.
\]
Let
\[
	E^\zeta~:\quad Y^2=X^3+\delta_\zeta X^2 -X.
\]
As above, $E$ will denote the elliptic curve 20A2.

\begin{lem}\label{lem:Fibpm1}
	Suppose $p \geq 7$ is a prime. 
	Suppose there exists an integer $k$
	satisfying the following conditions:
	\begin{enumerate}
		\item[(a)] The integer $l=2kp+1$ is prime, and $l \equiv \pm 1 \pmod{5}$.
		\item[(b)] The order of $\omega$ modulo $l$ is divisible by $p$; equivalently
				$\omega^{2k} \not \equiv 1 \pmod{l}$.
		\item[(c)] For all $\zeta \in A(p,k)$, we have
				\[
					a_l(E^\zeta)^2 \not \equiv a_l(E)^2 \pmod{p}.
				\]
	\end{enumerate}
	Then any solution to the equation~(\ref{eqn:FibPrimeIndex}) 
	must satisfy $n \equiv \pm 1 \pmod{p}$.
\end{lem}
\begin{proof}
	Suppose $p,~k$ satisfy the conditions of the Lemma, and that
	$(n,y,p)$ is a solution to equation~(\ref{eqn:FibPrimeIndex}). 
	Let $H_n$ and $E_n$
	be as above. Thus $H_n$ satisfies~(\ref{eqn:hn2}).

	We will prove first that $l \nmid y$. Suppose that $l \mid y$. Then
	$(\omega^n+\omega^{-n})/\sqrt{5}= F_n=y^p \equiv 0 \pmod{l}$ 
	and so $\omega^{4n} \equiv 1 \pmod{l}$.  From (b) 
	we deduce that $p \mid 4n$. However, $n$ is prime, and so $p=n$.
	This is impossible, since otherwise $F_p=y^p$ and clearly $1 < F_p < 2^p$.
	Hence $l \nmid y$.

	Next we will show that $y^{2p} \equiv 1 \pmod{l}$. Thus suppose that
	$y^{2p} \not \equiv 1 \pmod{l}$. Then there is some $\zeta \in A(p,k)$
	such that $y^{2p} \equiv \zeta \pmod{l}$. Further
	$\delta_\zeta \equiv \pm H_n \pmod{l}$. It follows that
	$a_l(E^\zeta)=\pm a_l(E_n)$.
	Applying Proposition~\ref{prop:FibLL} again, 
	we see that $a_l(E_n) \equiv a_l(E) \pmod{p}$.
	These congruences now contradict condition (c).

	We have finally proven that $y^{2p} \equiv 1 \pmod{l}$. By equation~(\ref{eqn:hn2})
	we see that $H_n \equiv \pm 1 \pmod{l}$. Since $n$ is odd (in fact an odd prime),
	and $\tau=-\omega^{-1}$, we get from the definition of $H_n$ that
	$\omega^{2n} \pm \omega^n -1 \equiv 0 \pmod{l}$. Solving this we find
	that $\omega^n \equiv \pm \omega^{\pm 1} \pmod{l}$. Thus
	\[
		\omega^{2(n+1)} \equiv 1 \pmod{l} \qquad
	\text{or}
		\qquad \omega^{2(n-1)} \equiv 1 \pmod{l}.
	\]
	However, condition (b) of the Lemma assures us that the order of $\omega$
	modulo $l$ is divisible by $p$. This immediately shows that $n \equiv \pm 1 \pmod{p}$
	as required.
\end{proof}

\begin{proof}[Proof of Proposition~\ref{prop:Fibnequivpm1}]
	We used a {\tt PARI/GP} program to check that for each prime
	$7 \leq p < 2 \times {10}^8$, there is some $k$ satisfying conditions (a),(b),(c)
	of Lemma~\ref{lem:Fibpm1}. This took approximately 41 hours on a $1.7$ GHz
	Pentium 4.
	This proves the Proposition. 
\end{proof}

\section{Level-lowering for Lucas - The Odd Index Case}
In this section we associate a Frey curve to solutions of 
~(\ref{eqn:LucPrimeIndex}) and apply level-lowering. Our
objective is to give the Lucas analogue of 
Propositions~\ref{prop:FibLL} and~\ref{prop:Fibnequivpm1}.

Suppose then that $(n,y,p)$ is a solution to~(\ref{eqn:LucPrimeIndex}),
and we associate to this solution the Frey curve
\begin{equation} \label{eqn:LucFrey}
E_{n}: \qquad Y^2 = X^3- 5 F_n X^2+ 5 X.
\end{equation}

Let $E$ be the elliptic curve $200B1$ in Cremona's tables~\cite{Cre}. 
This has the model
\begin{equation}\label{eqn:LucEll}
E: \quad Y^2=X^3+ X^2 -3X -2.
\end{equation}
\begin{prop}\label{prop:LucLL}
Suppose that $(n,y,p)$ is a solution to~(\ref{eqn:LucPrimeIndex}).
With notation as above, the Galois representations $\rho_p(E_n),~\rho_p(E)$
are isomorphic. Moreover, for any prime $l \neq 2,5$
\begin{itemize}
        \item[(i)] $a_l(E_n) \equiv a_l(E) \pmod{p}$ if $l \nmid y$.
        \item[(ii)] $l+1 \equiv \pm a_l(E) \pmod{p}$
         if $l \mid y$.
\end{itemize}
\end{prop}

\begin{prop}\label{prop:Lucnequivpm1}
        If $(n,y,p)$ is any solution to equation~(\ref{eqn:LucPrimeIndex})
        then $n \equiv \pm 1 \pmod{p}$.
\end{prop}

The proof of Proposition~\ref{prop:LucLL} is by no means as simple as
the proof of the corresponding Proposition for Fibonacci. 
However, given Proposition~\ref{prop:LucLL}, the proof of 
Proposition~\ref{prop:Lucnequivpm1} is a fairly trivial modification
of the proof of Proposition~\ref{prop:Fibnequivpm1} and we omit it.
The reader will notice that in Proposition~\ref{prop:Fibnequivpm1} 
(the Fibonacci case) we suppose that $p< 2 \times {10}^8$, but in the Lucas case
above there is no such assumption. This is because we know by a result of 
Peth\H{o} quoted in Section~\ref{sec:survey} that $p<13222$, which also
means that our program for the proof of Proposition~\ref{prop:Lucnequivpm1}
takes only a few seconds. Later on we will prove a much better bound for
$p$ in the Lucas case, but this will be dependent on the fact that 
$p \equiv 1 \pmod{n}$ and hence cannot be used for the proof of 
Proposition~\ref{prop:Lucnequivpm1}.

\subsection{Level-Lowering}
Let $E^1,\ldots,E^5$ be the elliptic curves 200A1,
200B1, 200C1, 200D1, 200E1 in Cremona's tables \cite{Cre}.
Note that $E^2$ is just our elliptic curve $E$ defined above.
We follow the notation of previous sections  with
regard to Galois representations and traces of Frobenius. 

\begin{lem}\label{lem:LucLL}
Suppose $(n,y,p)$ is a solution to equation~(\ref{eqn:LucPrimeIndex}). 
With notation as above, the Galois representations $\rho_p(E_n)$
is isomorphic to one of the Galois representations
$\rho_p(E^1),\ldots,\rho_p(E^5)$.
Moreover, if $\rho_p(E_n)$ is isomorphic to $\rho_p(E^i)$ then,
for any prime $l \neq 2,~5$,
\begin{itemize}
	\item[(i)] $a_l(E_n) \equiv a_l(E^i) \pmod{p}$ if $l \nmid y$. 
	\item[(ii)] $l+1 \equiv \pm a_l(E^i) \pmod{p}$
		 if $l \mid y$.
\end{itemize}
\end{lem}
\begin{proof}
	By the results of \cite[Sections 2,~3]{BS},
	$\rho_p(E_n)$ arises from a cuspidal newform of
	weight $2$,  level $200$, and trivial Nebentypus character. For this we need
	Lemma~\ref{lem:LucAux}. However using
	{\tt Magma} we find that the dimension of newforms of weight $2$ and level $200$
	is $5$.
	However, there are (up to isogeny) exactly $5$ elliptic curves of conductor $200$,
	and these are the curves $E^1,\ldots,E^5$ above.

	The rest of the
	lemma follows from \cite[Proposition 3]{KO}, and the fact that if
	$l \neq 2,5$ and $l \mid y$ then $l$ is a prime of multiplicative reduction
	for $E_n$ and so $a_l(E_n)= \pm 1$.
\end{proof}

\subsection{Eliminating Newforms}
Lemma~\ref{lem:LucLL} relates the Galois representation of $E_n$
to too many Galois representations. We now eliminate all but
one of them.

Suppose $l \neq 2,~5$ is a prime.  Define $d_l(E_n,E^i)=a_l(E_n)-a_l(E^i)$.
Let $M(l)$ be given by~(\ref{eqn:M(l)}). Recall that (Lemma~\ref{lem:M(l)}) the
residue class of $F_n$ modulo $l$, and hence the Frey curve $E_n$ modulo $l$,
depends only on the residue
class of $n$ modulo $M(l)$. We see that the following
definitions make sense: let
\[
\Tc_l(E^i) = \left\{m \in \Z/M(l): d_l(E_m,E^i)=0  \right\},
\]
\[
g_l(E^i)= {\rm lcm} \left\{ d_l(E_m,E^i): m \in \Z/M(l), \quad m \not \in \Tc_l (E^i) \right\},
\]
and
\[
h_l(E^i)= \left\{
\begin{array}{ll}
g_l(E^i) & \text{if $l \equiv \pm 2 \pmod{5}$}\\
{\rm lcm}(g_l(E^i),l+1-a_l(E^i),l+1+a_l(E^i))  & \text{if $l \equiv \pm 1 \pmod{5}$.}
\end{array}
\right.
\]

\begin{lem}\label{lem:eliminate}
Suppose that $l \neq 2,~5$ is a prime.
If $\rho_p(E_n)$ is isomorphic to $\rho_p(E^{(i)})$
then either the reduction of $n$ modulo $M(l)$ belongs to $\Tc_l(E^i)$
or else $p$ divides $h_l(E^i)$.
\end{lem}
\begin{proof}
Recall that by Lemma~\ref{lem:LucAux}, $y^{2p}=5F_n^2-4$. Thus
if $l \equiv \pm 2 \pmod{5}$ then $l$ does not divide $y$.
The Lemma now follows from Lemma~\ref{lem:LucLL}.
\end{proof}

Given two positive integers $M_1,~M_2$, and two sets
$T_1 \subset \Z/{M_1}$ and $T_2 \subset \Z/{M_2}$
we loosely define their \lq intersection\rq\ $T_1 \cap T_2$ to be the set of
all elements of $\Z/{{\rm lcm}(M_1,M_2)}$
whose reduction modulo $M_1$ and $M_2$ is respectively in $T_1$ and $T_2$.

We are now ready to prove Proposition~\ref{prop:LucLL}.
\begin{proof}[Proof of Proposition~\ref{prop:LucLL}]
	Suppose that $(n,y,p)$ is a solution to~(\ref{eqn:LucPrimeIndex}).
	Thus $p \geq 7$ and $n \equiv \pm 1 \pmod{6}$.
	We recall that the elliptic curves  $E$ and $E^2$ are one and the same. 
	Thus the Proposition follows from Lemma~\ref{lem:LucLL} if we can demonstrate
	that $\rho_p(E_n)$ cannot be isomorphic to the corresponding representation
	for $E^1,E^3,E^4,E^5$. Fix $i$ one of $1,3,4,5$.
	By the above Lemma, to show that the Galois representations of $E_n$
	and $E^i$ are not isomorphic it it is enough to produce a set
	of primes
	$S=\left\{l_1,\ldots,l_r \right\}$ all neither $2$ nor $5$ satisfying
	\begin{enumerate}
		\item For every $l \in S$ the integer $h_{l}(E^i)$ is not divisible
   			by any prime greater than $5$,
	\item $\left( \cap_{l \in S} \Tc_{l}(E^i) \right) \cap \Tc_0 = \emptyset$,
	\end{enumerate}
	where $\Tc_0 = \left\{ \overline{1},\overline{5} \right\} \subseteq \Z/6Z$.
	With the help of a short {\tt PARI/GP} program we find that
	we can take $S=\left\{ 3 \right\}$ to eliminate $E^1,E^3,E^5$ and
	$S=\left\{ 3,7,11,13,17,19,23 \right\}$ to eliminate $E^4$.

	We note in passing that the $j$-invariant of the curve $E^3$ is $55296/5$,
	and so the argument used in the proof of Lemma~\ref{lem:LucCong} also shows
	that the Galois representation $\rho_p(E^3)$ is not isomorphic to $\rho(E_n)$. This
	argument does not apply to the Galois representations of $E^1,E^4,E^5$ 
	as these have integral $j$-invariants.
\end{proof}

\section{Bounds for $n$ in terms of $p$}

Our objective in this section is to obtain bounds for $n$
in terms of $p$ for solutions to~(\ref{eqn:FibPrimeIndex})
and~(\ref{eqn:LucPrimeIndex}).
It follows from Baker's theory of linear
forms in logarithms (see for example the book of Shorey and Tijdeman \cite{ST}) that
the sizes of $n$ and $y$ are bounded in terms of $p$.
Unfortunately, these bounds are
huge, and there is no hope to complete the resolution of our
equations by proceeding in that way.
We however recall, by Lemma~\ref{lem:FibAux} (and~\ref{lem:LucAux}), 
that it is sufficient to obtain upper bounds
for the size of integer solutions to the equation $x^2+4=5y^{2p}$
(and one like it in the Lucas case).
As is explained below, this
equation easily reduces to a Thue equation, and we may apply
the results of Bugeaud and Gy\H ory \cite{BG} to get an upper bound for $x$
and $y$. However, it is of much interest to rework the proof
of Bugeaud and Gy\H ory in our particular context. On the one hand, our
particular equation has some nice properties not taken into
account in the general result of \cite{BG}, and, on the other hand,
there has been an important improvement, due to Matveev, in the
theory of linear forms in logarithms since \cite{BG} has appeared.
Altogether we actually compute a much better
upper bound, than the one obtained by applying directly the main
result of \cite{BG}. 

Before giving a precise statement of the main results of this section, 
we need a upper bound for the regulators of number fields. 
Several explicit upper bounds for regulators of a number field are available in
the literature; see for example \cite{Len} and \cite{Siegel}.
We have however found it best to use a result of Landau.

\begin{lem}~\label{lem:Landau}
	Let $\K$ be a number field
        with degree $d=r_1+2r_2$
	where $r_1$ and $r_2$ are numbers of real and complex embeddings.
	Denote the discriminant by $D_\K$ and the regulator by $R_\K$,
	and the number of roots of unity in $\K$ by $w$. Suppose, moreover,
	that $L$ is a real number such that $D_\K \leq L$. 
	Let
	\[
	     a = 2^{-r_2} \, \pi^{-d/2} \, \sqrt{L}.
	\]
	Define the function $f_{\K}(L,s)$ by
	\[
	      f_{\K}(L,s)=2^{-r_1} \, w \,
	       a^s \, \bigl( \Gamma (s/2) \bigr)^{r_1} \, \bigl( \Gamma(s) \bigr)^{r_2}
	       s^{d+1} \, (s-1)^{1-d},
	\]
	and let $C_{\K}(L)={\rm min} \left\{f_{K}(L,2-t/1000)~: t=0,1,\ldots,999 \right\}$.
	Then $R_{\K} < C_\K(L)$.
\end{lem}
\begin{proof}
		Landau \cite{Landau} proved that $R_{\K} < f_{\K} (D_{\K},s)$ 
		for all $s>1$. It is thus
		clear that $R_{\K} < C_{\K}(L)$.

		Perhaps a comment is in order. For a complicated number 
		field of high degree it is difficult to calculate the discriminant 
		$D_{\K}$ exactly, though it is easy to give an upper bound $L$ for
		its size. It is also difficult to
		minimise the function $f_{\K}(L,s)$ analytically, but we have 
		found that the above gives
		an accurate enough result, which is easy to 
		calculate on a computer.
\end{proof}

We are now ready to state our upper bound for $n$ in terms of $p$ for the
Fibonacci and Lucas cases.
\begin{prop}\label{prop:Bugeaud1}
        Suppose $p \geq 7$ is prime.
	Let $\alpha$ be any root of the polynomial
	\begin{equation}\label{eqn:P}
	        P(X) := \sum_{k=0}^p \, 
	        (-4)^{[(p-k)/2]} \, \binom{p}{k} \,  X^k,
	\end{equation}
	and let $\K=\Q(\alpha)$. Let $C_{\K}(\cdot)$ be as in Lemma~\ref{lem:Landau}
	and 
        \[ 
               \Theta =  3.9 \cdot 30^{p+3} \, p^{13/2} \, (p-1)^{p+1} \, \bigl((p-1)!\bigr)^2 \,
               (3p + 2) \, \bigl(1 + \log (p(p-1)) \bigr) \, C_{\K}(10^{p-1} p^p).
        \]
	If $(n,y,p)$ satisfies the equation and conditions~(\ref{eqn:FibPrimeIndex}) 
	then $n < 2.5 p \Theta \log{\Theta}$.
\end{prop} 

\begin{prop}\label{prop:Bugeaud2}
        Suppose $p \geq 7$ is prime.
	Denote by $\root p\of{\omega}$ the real $p$-th root of $\omega$ (where $\omega$
	is given by~(\ref{eqn:lmu})) and set
	$\K = \Q(\sqrt{5}, \root p\of{\omega})$, and let $C_\K(\cdot)$ 
	be as in Lemma~\ref{lem:Landau}.
	Let
	\[
		\Theta =  67 \cdot 30^{p+5} (p-1)^{p+2} p^3 \, (p+2)^{5.5} \, (p!)^2 \,
		\bigl(1 + \log (2p(p-1)) \bigr)  \, C_\K(5^p p^{2p}).
	\]
	If $(n,y,p)$ satisfies the equation and conditions~(\ref{eqn:LucPrimeIndex}) 
	then $n < 2.5 p \Theta \log \Theta$.
\end{prop}

\subsection{Preliminaries}
We first need a lower bound for linear forms in logarithms, due to 
Matveev. Let $\LL$ be a number field of degree $D$, let $\al_1, \ldots, \al_n$
be non-zero elements of $\LL$ and $b_1, \ldots, b_n$ be rational integers. Set
\[
	B = \max\{|b_1|, \ldots, |b_n|\},
\]
and
\[
	\La = \al_1^{b_1} \ldots \al_n^{b_n} - 1.
\]
Let $\h$ denote the absolute logarithmic height and let
$A_1, \ldots, A_n$ be real numbers with
\[
	A_j \ge \h'(\al_j) := 
	\max\{D \, \h(\al_j), |\log \al_j|, 0.16\}, \quad 1 \le j \le n.
\]
We call $\h'$ the modified height (with respect to the field $\LL$).
With this notation, the main result of Matveev \cite{Mat}
implies the following estimate.

\begin{theorem}\label{thm:M} 
	Assume that $\La$ is non-zero. We then have
	\[
		\log |\La| > - 3 \cdot 30^{n+4} \, (n+1)^{5.5} \, 
		D^2 \, (1 + \log D) \, (1 + \log nB) \, A_1 \ldots A_n.
	\]
	Furthermore, if $\LL$ is real, we have
	\[
		\log |\La| > - 1.4 \cdot 
		30^{n+3} \, n^{4.5} \, D^2 \, (1 + \log D) \, (1 + \log B) \, A_1 \ldots A_n.
	\]
\end{theorem}
\begin{proof}
	Denote by $\log$ the principal determination of the logarithm. If
	$|\Lambda| < 1/3$, then there exists an integer $b_0$, with 
	$|b_0| \le n\, B$, such that 
	\[
	\Omega := |b_0 \log (-1) + b_1 \log \al_1 + \ldots + b_n \log \al_n| 
	\]
	satisfies $|\Lambda| \ge \Omega / 2$. Noticing that $\h'(-1) = \pi$,
	and that $b_0 = 0$ if $\LL$ is real, we
	deduce our lower bounds from Corollary 2.3 of Matveev \cite{Mat}.
\end{proof}

We also need some precise results from algebraic number theory.
In the rest of this Section, let $\K$ denote a number field of degree 
$d=r_1 + 2 r_2$ and unit rank $r= r_1 + r_2 - 1$ with $r > 0$.
Let $R_{\K}$ and $D_{\K}$ be its
regulator and discriminant, respectively. Let $w$ denote the 
number of roots of unity in $\K$. Observe that $w=2$ if $r_1 > 0$.

\begin{lem}\label{lem:B1}
	For every algebraic integer $\eta$ which generates
	$\K$ we have 
	\[
		d \, \h(\eta) \ge \frac{\log |D_{\K}| - d \, \log d}{2 (d-1)}. 
	\]	
\end{lem}
\begin{proof}
	As in Mignotte \cite{Mi3}, it follows from the Hadamard inequality that
	\[ 
	|D_{\K}|
	\le {\rm Discr}(1, \eta, \ldots, \eta^{d-1})^2 \le d^d \, \M(\eta)^{2(d-1)}, 
	\]
	where $\M(\eta)$ is the Mahler measure of $\eta$.
	Since $d \, \log \M(\eta) = \h(\eta)$, the lemma is proved. 
\end{proof}

In the course of our proof, we use fundamental systems of units in $\K$ with
specific properties.

\begin{lem}\label{lem:BG} 
	There exists in $\K$ a fundamental system 
	$\{\varepsilon _1, \ldots, \varepsilon _r\}$ of units such that
	\[
	\prod_{i=1}^{r} \, \h(\varepsilon _i) \le 
	2^{1 - r} \, ( r! )^2 \, d^{-r} \, R_{\K},
	\]
	and the absolute values of the entries of 
	the inverse matrix of $(\log |\varepsilon _i|_{v_j})_{i, j = 1, \ldots, r}$
	do not exceed $(r!)^2 \, 2^{-r} \, (\log (3d))^3$.
\end{lem}
\begin{proof}
	This is Lemma 1 of \cite{BGunits}
	combined with a result of Voutier \cite{Vou} (see \cite{BG})
	giving a lower bound for the height of any non-zero algebraic
	number which is not a root of unity. 
\end{proof}

Furthermore, we need sharp bounds for discriminants of number
fields in a relative extension.

\begin{lem}\label{lem:DiscRel} 
	Let $\K_1$ and $\K_2$ be number fields with $\K_1 \subseteq \K_2$,
	and denote the discriminant of the extension
	$\K_2/\K_1$ by $D_{\K_2 / \K_1}$. Then
	\[
	|D_{\K_2}| = |D_{\K_1}|^{[\K_2 : \K_1]} \, 
	|{\rm N}_{\K_1 / \Q} (D_{\K_2 / \K_1})|.
	\]
\end{lem}
\begin{proof}
	This is Proposition 4.9 of \cite{Nark}.
\end{proof}

\subsection{Proof of Proposition~\ref{prop:Bugeaud1}}

We now turn our attention to the proof of Proposition~\ref{prop:Bugeaud1}
and so to equation~(\ref{eqn:FibPrimeIndex}). 
Lemma~\ref{lem:FibAux} reduces the problem to solving the 
superelliptic equation $x^2 + 4 = 5 y^{2p}$. Factorising the left-hand side
over $\Z[i]$, we deduce the existence of integers $a$ and $b$ with 
$a^2 + b^2 = y^2$ and
\begin{equation}\label{eqn:B3}
	\pm 4i = (2 + i) (a + ib)^p - (2 -i) (a - ib)^p .
\end{equation}
Dividing by 2i, we get
\begin{equation*}\label{eqn:B4}
	\begin{split}
	\pm 2 = & 2 \, \sum_{k=0}^{[p/2]} \, 
	\binom{p}{2k} \, a^{2k} \, (-1)^{(p-2k-1)/2} \, b^{p -2k} \\
	& + 
	\sum_{k=0}^{[p/2]} \, \binom{p}{2k+1}  \, a^{2k+1} 
	\, (-1)^{(p-2k-1)/2} \, b^{p-2k-1}.   
	\end{split}
\end{equation*}  
We infer  
that $a$ is even. Consequently, $(b, a/2)$ is an 
integer solution of the Thue equation
\begin{equation}\label{eqn:B4a}
	\sum_{k=0}^p \, (-4)^{[(p-k)/2]} \, \binom{p}{k} \,  X^k \, Y^{p-k} = \pm 1.   
\end{equation} 
To bound the size of the solutions of~(\ref{eqn:B4a})
we follow the general
scheme of \cite{BG}, which was also used in \cite{BMRS}.
Let $P(X)$ and $\alpha$ and $\K$ be as in Proposition~\ref{prop:Bugeaud1}; 
we note that $P(X)$ is 
the polynomial naturally associated to the Thue equation~(\ref{eqn:B4a}).
We first need information on the number field $\K$ and its
Galois closure. We would like to thank Mr. Julien Haristoy
for his help in proving the following Lemma.

\begin{lem}\label{lem:FibDisc}
	The field $\K = \Q(\alpha)$ is totally real and 
	its Galois closure $\LL$ has degree $p(p-1)$ over $\Q$.
	Furthermore, the discriminant of $\K$ divides $10^{p-1} p^p$.
\end{lem}
\begin{proof}
	Observe that any root of the polynomial
	\[
		Q(X) := \frac{1}{2i} \cdot
		\bigl((2 + i) (X + i)^p - (2 - i) (X-i)^p\bigr)
		= (-1)^{(p-1)/2} (X/2)^p P(2/X).
	\]
	satisfies $|X + i| = |X-i|$, and so must
	be real. Hence, $\LL$ is a totally real field. 
	Furthermore, $\LL(i) / \Q(i)$ is a Kummer extension got by
	adjoining the $p$-th roots of unity and the $p$-th roots of
	$(2+i)/(2-i)$. Hence, this extension has degree $p(p-1)$, and
	this is the same for $\LL / \Q$.
	
	Observe now that $\K(i)$ is generated over $\Q(i)$ by any 
	root of either of the following two monic polynomials
	with coefficients in $\Z[i]$, namely $Y^p - (2+i) (2-i)^{p-1}$ and
        $Y^p-  (2-i) (2+i)^{p-1}$. Since the discriminant $D_1$ (viewed as
        an algebraic integer in $\Z[i]$ and not as an ideal) of the extension
        $\K(i)/\Q(i)$ divides the discriminant of each of these polynomials,
        $D_1$ divides $p^p 5^{p-1} (2-i)^{(p-1)(p-2)}$ 
        and $p^p 5^{p-1} (2+i)^{(p-1)(p-2)}$.
        However, $2+i$ and $2-i$ are relatively prime, thus $D_1$ divides 
        $5^{p-1} p^p$. Furthermore, estimating the 
        discriminant of $\K(i) / \Q$ in two different ways
        thanks to Lemma \ref{lem:DiscRel} gives
       	\begin{equation}\label{eqn:discRel} 
        |D_{\K(i)}| = 4^p D_1^2 = |D_{\K}|^2 \cdot  |{\rm N}_{\K/\Q} 
        (D_{\K(i)/\K})|.  
        \end{equation} 
        Consequently, $|D_{\K}|$ divides $5^{p-1} (2p)^p$.
        We now refine this estimate by showing that $4$ divides 
        $|{\rm N}_{\K/\Q} (D_{\K(i)/\K})|$.
        
        Suppose that the decomposition of the ideal $2 \cdot \OO_{\K}$ in $\K / \Q$ is given by
        \[
        2 \cdot \OO_{\K} = \pp_1^{e_1} \ldots \pp_s^{e_s}.
        \]
        At least one of the $e_i$ is odd, since otherwise $2$ would 
        divide $\sum_{i=1}^s e_i f_i = p$. Thus, there is (at least)
        one prime $\pp$ in $\OO_{\K}$ lying above $2$ whose ramification
	index $e$ is odd: this prime must ramify in $\K(i) / \K$, since
	$2i=(1+i)^2$ in $\K(i)$. 
	Thus $\pp$ divides $D_{\K(i)/\K}$ and so $|{\rm N}_{\K/\Q} (D_{\K(i)/\K})|$
	is divisible by $2$. However, by~(\ref{eqn:discRel}), we
	know that $|{\rm N}_{\K/\Q} (D_{\K(i)/\K})|$ is a square and
	so must be divisible by $4$.
	
\end{proof}

\noindent {\bf Remark:} Based on computations for small $p$, 
it seems very likely that  $10^{p-1} p^p$ is the exact
value of $|D_{\K}|$ for most $p$.

Since we introduce many changes in the proof of \cite{BG}, we
give a complete proof, rather than only quoting \cite{BG}.

Let $\alpha_1, \ldots, \alpha_p$ be the roots of $P(X)$ and
let $(X, Y)$ be a solution of~(\ref{eqn:B4a}). Without any loss of generality,
we assume that $\al = \al_1$ and
$|X - \alpha_1 Y | = \min_{1 \le j \le p} |X - \al_j Y|$.
We will make repeated use of the fact that $|\al_1|, \ldots 
|\al_p|$ are not greater than $4^p$, neither 
smaller than $4^{-p}$ (since $4^p - 1$ is an
upper bound for the absolute values of the coefficients of $P(X)$).
Assuming that $Y$ is large enough, namely that
\begin{equation}\label{eqn:BoY}
	\log |Y| \ge (30 p)^p,
\end{equation}
we get $|Y| \ge 2 \, \min_{2 \le j \le p} \, \{|\al_1 - \al_j|^{-1}\}$ and
\begin{equation}\label{eqn:B5}
	|X - \al_1 Y| \le 2^{p-1} \, \prod_{2 \le j \le p} \, 
	|\al_1 - \al_j| \,|Y|^{-p+1} \le 2^{2 p^2} \, |Y|^{-p+1},  
\end{equation}
since $|X - \al_j Y| \ge |\al_1 - \al_j| \cdot |Y| / 2$ if
$|X - \al_1 Y| \le |\al_1 - \al_j| \cdot |Y| / 2$,
for any $j = 2, \ldots, p$.

{} From the `Siegel identity'
\[
	(X - \alpha_1 Y) (\alpha_2 - \alpha_3) + (X - \alpha_2 Y) (\alpha_3 - \alpha_1)
	+ (X - \alpha_3 Y) (\alpha_1 - \alpha_2) = 0
\]
we have
\[
	\La :=  \frac{\alpha_2 - \alpha_3}{\alpha_3 - \alpha_1} 
	\cdot  \frac{X - \alpha_1 Y}{X - \alpha_2 Y}  = 
	\frac{X - \alpha_3 Y}{X - \alpha_2 Y} \cdot
	\frac{\alpha_2 - \alpha_1}{\alpha_3 - \alpha_1} -1.
\]
Observe that the unit rank of $\K$ is $p-1$, since $\K$ is totally real.
Let $\eps_{1, 1}, \ldots, \eps_{1,p-1}$ be a fundamental system of units
in $\K := \Q(\alpha_1)$ given by Lemma~\ref{lem:BG}, hence, satisfying
\begin{equation}\label{eqn:B6}
	\prod_{1 \le i \le p-1} \, \h(\eps_{1,i}) \le 
	\frac{\bigl( (p-1)! \bigr)^2}{2^{p-2} p^{p-1}} \, R_{\K},  
\end{equation}
where $R_{\K}$ denotes the regulator of the field $\K$.
For $j=2, 3$, denote by 
$\eps_{2, 1}, \ldots, \eps_{2,p-1}$ and $\eps_{3, 1}, \ldots, \eps_{3,p-1}$
the conjugates of $\eps_{1, 1}, \ldots, \eps_{1,p-1}$
in $\Q(\alpha_2)$ and $\Q(\alpha_3)$, respectively.
They all belong to the Galois closure $\LL$ of $\K$.

The polynomial $P(X)$ is monic and the left-hand side of~(\ref{eqn:B4a}) 
is a unit, thus $X - \alpha_1 Y$ is a unit. 
This simple observation 
appears to be crucial, since, roughly speaking, it allows us to gain
a factor of size around $p^p R_{\K}$ (compare with the proofs in \cite{BG}
and in \cite{BMRS}).

Since the only roots of unity in $\K$ are $\pm 1$, 
there exist integers $b_1, \ldots, b_{p-1}$ such that
$X - \al_1 Y = \pm \eps_{1, 1}^{b_1} \ldots \eps_{1, p-1}^{b_{p-1}}$, thus
we have
\[
	\La = \pm \biggl( \frac{\eps_{3,1}}{\eps_{2,1} } \biggr)^{b_1}
	\ldots  \biggl( \frac{\eps_{3,p-1}}{\eps_{2,p-1} } \biggr)^{b_{p-1}} \,
	\frac{\alpha_2 - \alpha_1}{\alpha_3 - \alpha_1} - 1.
\]
As in \cite[6.12]{BG},
we infer from Lemma~\ref{lem:BG} that
\begin{equation}\label{eqn:B6b}
\begin{split}
	B := \max\{|b_1|, \ldots, |b_{p-1}|\} & 
	\le 2^{2-p} \, p \, (p!)^2 \,  (\log (3p))^3 \,
	\h(X - \al_1 Y) \\
	& \le p^{2(p+1)} \, \log |Y|,   
\end{split}
\end{equation}
by~(\ref{eqn:B5}).

Further, we notice that
\[
	\h \biggl( \frac{\alpha_2 - \alpha_1}{\alpha_3 - \alpha_1} \biggr)
	= \h \biggl( \frac{\alpha_2/2 - \alpha_1/2}
	{\alpha_3/2 - \alpha_1/2} \biggr)
	\le 4 \h(\alpha_1/2) + \log 4 \le \frac{6p+4}{p} \log 2,
\]
since we have (here and below, $M(\cdot)$ denotes the Mahler measure 
and ${\rm H} (\cdot )$ stands for the na\"\i ve height)
\[
	\h(\alpha_1/2) \le \frac{\log M(Q)}{p} 
	\le \frac{\log \bigr( \sqrt{p+1} {\rm H}(Q) \bigr)}{p} 
	\le \frac{\log \bigr( 2 \sqrt{p+1} \binom{p}{[p/2]} \bigr)}{p}
	\le \frac{p+1}{p} \, \log 2.
\]
Hence, with the modified height $\h'$ related to the field $\LL$, we have
\[
	\h' \biggl( \frac{\alpha_2 - \alpha_1}{\alpha_3 - \alpha_1} \biggr)
	\le 2(3p+2)(p-1) \log 2.
\]
We may assume from Lemma~\ref{lem:FibDisc} that the absolute value of 
the discriminant of $\K$ is $10^{p-1} p^p$, since the upper bound for
$n$ we aim to prove is an increasing function of $|D_{\K}|$.
For $i=1, \ldots, p-1$, we have 
$\h(\eps_{1,i}) = \h(\eps_{2,i})= \h(\eps_{3,i})$ and, by Lemma~\ref{lem:B1},
the height of the real algebraic integer
$\eps_{1,i}$ satisfies $\h(\eps_{1,i}) \ge (\log 10)/2$. Thus, we get
\[
	\h' \biggl( \frac{\eps_{2,i}}{\eps_{3,i} } \biggr) 
	\le 2 p (p-1) \h(\eps_{1,i}).
\]
Consequently, using Theorem~\ref{thm:M} in the real case
with $n= p$ and $D = p(p-1)$, we get
\begin{equation}\label{eqn:B7}
	\begin{split}
		\log |\La| & > - 1.4 \cdot 30^{p+3} p^{7/2} \bigl(p (p-1) \bigr)^{p+2} \, 
		(3p+ 2) \bigl( 1 + \log (p(p-1)) \bigr) \\ 
		& \times (1 + \log B) 
		\, (2 \log 2) \, 2^{p-1} \, \prod_{1 \le i \le p-1} \h(\eps_{1,i}). 
	\end{split}
\end{equation}
Then,~(\ref{eqn:B6}) gives us that
\begin{equation}\label{eqn:B8}
	\begin{split}
		\log |\La| & > - 3.9 \cdot 30^{p+3} \, 
		p^{13/2} \, (p-1)^{p+2} \, (3p+2) \, \bigl((p-1)!\bigr)^2 \\
		& \times \bigl( 1 + \log (p(p-1)) \bigr) \, (1 + \log B) \, R_{\K}. 
	\end{split}
\end{equation}
Furthermore, it follows from~(\ref{eqn:B5}) that
\begin{equation}\label{eqn:B9}
	\log |\La| < 5 p^2 - (p-1) \log |Y|. 
\end{equation}
By~(\ref{eqn:B6b}), we have the upper bound
\begin{equation}\label{eqn:B10}
	(1 + \log B) < 3 p^2 + \log \log |Y|.  
\end{equation}
Finally, we observe that if $F_n$ is a $p$-th power for some odd $n$, 
then there are integers $X$ and $Y$ such that $(X, Y)$ is a
solution of the Thue equation~(\ref{eqn:B4a}) and 
$F^{2/p}_n = 4 X^2 + Y^2$. Since $|X| \le 1 + 4^p |Y|$ and
$F_n \ge 0.4 \cdot 1.6^n$ (for $n \ge 7$), we derive from~(\ref{eqn:BoY})
that $n < 2.2 p \log |Y|$, 
It then follows from~(\ref{eqn:B8}), (\ref{eqn:B9}), and (\ref{eqn:B10}), 
together with Lemmas~\ref{lem:Landau} and~(\ref{lem:FibDisc}) that
\[
	n < 2.5 p \Theta \log \Theta,
\]
with
\[
	\Theta =  3.9 \cdot 30^{p+3} \, 
	p^{13/2} \, (p-1)^{p+1} \, (3p + 2) \, \bigl((p-1)!\bigr)^2 \, 
	\bigl(1 + \log (p(p-1)) \bigr) \, C_{\K}(10^{p-1} p^p). 
\]
This proves Proposition~\ref{prop:Bugeaud1}.

\subsection{Proof of Proposition~\ref{prop:Bugeaud2}}
Suppose that $(n,y,p)$ is a solution to the equation~(\ref{eqn:LucPrimeIndex}).
In particular, we know  
\[
	y^p=L_n = \omega^n + \tau^n,
\]
where we recall that $\omega = (1 + \sqrt{5})/2$ and $\tau$ is
conjugate of $\omega$. 
We also know by Proposition~\ref{prop:Lucnequivpm1} that $n$ is
congruent to $\pm 1$ modulo $p$. This means that there exists
an integer $\nu$ such that
\[
	y^p - \omega^{\pm 1} \, (\omega^{\nu})^p = - \tau^n.
\]
Thus, we are left with an equation of Thue type, namely
\begin{equation}\label{eqn:L10}
	X^p - \omega^{\pm 1} \, Y^p = \text{unit in $\Q(\sqrt{5})$}.  
\end{equation}
We only deal with the $+$ case, since the $-$ case is entirely similar.

As in the statement of Proposition~\ref{prop:Bugeaud2},
denote by $\root p\of{\omega}$ the real $p$-th root of $\omega$ and set
$\K = \Q(\sqrt{5}, \root p\of{\omega})$. Let $\zeta$ be a 
primitive $p$-th root of unity.

\begin{lem}\label{lem:L3} 
	The field $\K$ has degree $2p$ and we have $r_1 = 2$,
	$r_2 = p-1$ and $r=p$. The absolute value of the discriminant
	of $\K$ is at most equal to $5^p p^{2p}$. Its non-trivial subfields
	are $\Q(\sqrt{5})$ and $\Q(\root p\of{\omega} - (\root p\of{\omega})^{-1})$,
	whose discriminant is, in absolute value,
	at most equal to $5^{(p-1)/2} p^p$.
	Furthermore, the Galois closure $\LL$ of $\K$ is the field $\K(\zeta)$,
	of degree $2p (p-1)$.
\end{lem}
\begin{proof}
	We observe that the minimal defining polynomial of
	$\root p\of{\omega}$ over $\Z$ is $R(X) := X^{2p} - X^p - 1$, thus
	we have
	\[
	|D_{\K}| \le |{\rm N}_{\K / \Q} (R'(\root p\of{\omega})) |
	= |{\rm N}_{\K / \Q} ( p \sqrt{5} (\root p\of{\omega})^{p-1}) |
	= 5^p \, p^{2p}.
	\]
	The fact that $\K$ has only two non-trivial subfields, one of degree
	$2$, another of degree $p$, is clear.
	Furthermore, since $\K$ is obtained from 
	$\Q(\root p\of{\omega} - (\root p\of{\omega})^{-1})$ by adjoining
	$\sqrt{5}$, we get from 
	Lemma \ref{lem:DiscRel} that the absolute value of the discriminant of 
	$\Q(\root p\of{\omega} - (\root p\of{\omega})^{-1})$ is not greater
	than $5^{(p-1)/2} p^p$.
	Since the roots of the polynomial $R(X)$ are the algebraic numbers
	$\root p\of{\omega}, \zeta \root p\of{\omega}, \ldots,
	\zeta^{p-1} \, \root p\of{\omega}, \root p\of{\tau},
	\zeta \root p\of{\tau}, \ldots, \zeta^{p-1} \,
	\root p\of{\tau}$, the Galois closure of $\K$ is 
	the field $\K(\zeta)$.
\end{proof}

Let $\eps_{1,1}, \ldots, \eps_{1,p}$
be a fundamental system of units in $\K$ given by Lemma~\ref{lem:BG}.
There exist integers $b_1, \ldots, b_p$ such that
\[
	X - \root p\of{\omega} Y = \pm \eps_{1,1}^{b_1} \ldots \eps_{1,p}^{b_p}.
\]
Keep in mind that we are only interested in solutions $(X, Y)$ 
of~(\ref{eqn:L10})
with $X$ integer and $Y$ algebraic integer in the field $\Q(\sqrt{5})$.
Thus, $X/Y$ is real, $|X - \root p\of{\omega} Y|$ is small, and 
$|X - \zeta^j \root p\of{\omega} Y|$ is quite large for 
$j= 1, \ldots, p-1$ (consider the imaginary part). More
precisely, for $Y > 2$, we get
\begin{equation}\label{eqn:L11}
	|X - \root p\of{\omega} Y| \le p^p \, Y^{-p+1}.  
\end{equation}
Furthermore, setting $B = \max\{|b_1|, \ldots, |b_p|\}$,
Lemma~\ref{lem:BG} yields that
\begin{equation}\label{eqn:L12}
	\begin{split}
		B & \le 2^{1-p} \, p (p!)^2 (\log 6p)^3 \, 
			\h (X - \root p\of{\omega} Y) \\
		& \le p^{2(p+1)} \, \log Y,   
	\end{split}
\end{equation}
by our assumptions on $X$ and $Y$.

Recall that $\zeta$ is a primitive $p$-th root of unity.
We introduce the quantity
\begin{equation}\label{eqn:L13}
	\La :=  \frac{\zeta - \zeta^2}{\zeta^2 - 1} 
	\cdot  \frac{X - \root p\of{\omega} Y}{X - \zeta \root p\of{\omega} Y}  = 
	\frac{X - \zeta^2 \root p\of{\omega} Y}{X - \zeta \root p\of{\omega} Y} \cdot
	\frac{\zeta - 1}{\zeta^2 - 1} -1, 
\end{equation}
hence, the linear form in logarithms
\[
	\La = \biggl( \frac{\eps_{3,1}}{\eps_{2,1} } \biggr)^{b_1}
	\ldots  \biggl( \frac{\eps_{3,p}}{\eps_{2,p} } \biggr)^{b_{p}} \,
	\frac{\zeta - 1}{\zeta^2 - 1} - 1.
\]
Let $\h'$ denote the modified height related to the field $\LL$. We have 
\[
	\h' \biggl( \frac{\zeta - 1}{\zeta^2 - 1} \biggr) \le 2p(p-1) \log 4,
\]
and $\h'(\eps_{1, i}) = \h(\eps_{1, i})$. To check this, we observe that 
any algebraic unit in $\K$ generates one of the subfields of $\K$,
and we apply Lemma~\ref{lem:L3}
(we may assume that the absolute value of 
the discriminant of $\K$ is $5^p p^{2p}$, since the upper bound for
$n$ we aim to prove is an increasing function of $|D_{\K}|$).
Using Theorem~\ref{thm:M} in the complex case with $n=p+1$ and $D=2p(p-1)$, we get
\begin{equation}\label{eqn:L17}
	\begin{split}
		\log |\La| 
		& > - 3 \cdot 30^{p+5} (p+2)^{5.5} \bigl(2 p (p-1) \bigr)^{p+3} \,
		\bigl(1 + \log (2p(p-1)) \bigr) \\
		& \times (1 + \log (p+1)B) 
		\, (\log 4) \, 2^{p} \, \prod_{1 \le i \le p} \h(\eps_{1,i}). 
	\end{split}
\end{equation}

By~(\ref{eqn:L17}) and Lemma~\ref{lem:BG}, we get
\begin{equation}\label{eqn:L19}
	\begin{split}
		\log | \La| & > - 3 \cdot 30^{p+5} (p+2)^{5.5} \bigl(2 p (p-1) \bigr)^{p+3} \,
		\bigl(1 + \log (2p(p-1)) \bigr) \\ 
		& \times \bigl(1 + \log ((p+1)B) \bigr)
		\,(\log 4) \, 2^{-p+1} \, p^{-p} \, (p!)^2 \, R_{\K}.
	\end{split}
\end{equation}
Furthermore, it follows from~(\ref{eqn:L11}) and~(\ref{eqn:L13}) that
\begin{equation}\label{eqn:L20}
	\log |\La| < 5 p^2 - (p-1) \log |Y|. 
\end{equation}
Observe now that if $L_n = y^p$ for some $n$, then equation~(\ref{eqn:L10})
has a solution $(X, Y)$ with $Y = \omega^{(n \pm 1)/p}$, and we get
that $-1 < L_n - \omega^{\pm 1} y^p < 0$, thus $n < 2.2 p \log Y$.
It then follows from~(\ref{eqn:L12}),~(\ref{eqn:L17})--(\ref{eqn:L20}),
together with Lemma~\ref{lem:Landau} that
\[
	n < 2.5 p \Theta \log \Theta,
\]
with
\[
	\Theta =  67 \cdot 30^{p+5} (p-1)^{p+2} p^3 \, (p+2)^{5.5} \, (p!)^2 \,
	\bigl(1 + \log (2p(p-1)) \bigr)  \, C_{\K}(5^p p^{2p}). 
\]
This completes the proof of Proposition~\ref{prop:Bugeaud2}.

\section{The Sieve}\label{sec:sieve}
In this section we use Propositions~\ref{prop:FibLL} and~\ref{prop:Bugeaud1}
(the Fibonacci case) and Propositions~\ref{prop:LucLL} and~\ref{prop:Bugeaud2}
(the Lucas case) together with a substantial computation to prove the following.
\begin{prop}\label{prop:Fib733}
	If $(n,y,p)$ satisfies the equation and conditions~(\ref{eqn:FibPrimeIndex})
	then
	\[
	p > 733, \qquad n \geq 1.033 \times 10^{8733}, \qquad \log{y} > 10^{8000} .
	\]
\end{prop}
\begin{prop}\label{prop:Luc283}
	If $(n,y,p)$ satisfies the equation and conditions~(\ref{eqn:LucPrimeIndex})
	then
	\[
	p > 283, \qquad n \geq 4.938 \times 10^{3383}, \qquad  \log{y} > 10^{3000}.
	\]
\end{prop}
We will focus on the Fibonacci case; the Lucas case is entirely similar. Throughout
this section we will follow the notation of Section~\ref{sec:FibLL}. In particular,
$H_n,~E_n$ and $E$ are given respectively by~(\ref{eqn:Hn}),~(\ref{eqn:FibFrey}), 
and~(\ref{eqn:FibEll}).

\begin{lem}\label{lem:Kl}
Suppose $l \equiv \pm 1 \pmod{5}$ is prime and let
\[
K(l)= \rm{lcm}(l-1,6).
\]
The trace of Frobenius $a_l(E_n)$ depends only on 
the residue class of $n$ modulo $K(l)$. 
\end{lem}
\begin{proof}
By Lemma~\ref{lem:M(l)}, the residue class of $L_n$ modulo $l$ depends only on
the residue class of $n$ modulo $l-1$. From the definition of $H_n$ in~(\ref{eqn:Hn})
we see that $H_n$ modulo $l$ depends only on the residue class of $n$ modulo $K(l)$.
The Lemma follows at once from the fact that the Frey curve $E_n$ depends only on $H_n$.
\end{proof}
Suppose $l \equiv \pm 1 \pmod{5}$; we see by Lemma~\ref{lem:Kl} 
that for $n \in \Z/K(l)$ it makes sense to talk 
of $a_l(E_n)$.  Suppose $q \geq 5$ is a fixed prime. 
Define $\N(l,q)$ to be the subset of all $n \in (\Z/K(l))^*$
such that 
\begin{itemize}
	\item either $H_n^2+4 \not \equiv 0 \pmod{l}$, 
	and the integer $a_l(E_n) - a_l(E)$ is divisible by some prime $p > q$. 
	\item or $H_n^2+4 \equiv 0 \pmod{l}$ and  
	one of the two integers  $l+1 \pm a_l(E)$ is divisible by some prime
	$p > q$.
\end{itemize}

\begin{lem}\label{lem:Nl}
Suppose that $q \geq 5$ is prime.
Suppose $l$ satisfies
\begin{equation}\label{eqn:condl}
 l \equiv \pm 1 \pmod{5} \text{ is prime and every prime factor of $l-1$ is $< 25000$}.
 \end{equation}
If $(n,p,y)$ satisfies the equation~(\ref{eqn:FibPrimeIndex})
and $p > q$ then the reduction of $n$ modulo $K(l)$ belongs to $\N(l,q)$. 
\end{lem}
\begin{proof}
First observe, since $n$ satisfies~(\ref{eqn:FibPrimeIndex}),
that $n$ is prime and $n \geq 25000$. However, every prime
divisor of $l-1$ is $< 25000$ and the same must be true
of $K(l)={\rm lcm}(l-1,6)$. Thus the reduction of $n$
modulo $K(l)$ certainly belongs to $(\Z/K(l))^*$.

Next we recall (Lemma~\ref{lem:1mod4}) that $H_n^2+4=5y^{2p}$
and so $l|y$ if and only if $H_n^2+4 \equiv 0 \pmod{l}$.
The Lemma now immediately follows from Proposition~\ref{prop:FibLL}.
\end{proof}

Given two positive integers $M_1,~M_2$, and two sets 
$T_1 \subset \Z/{M_1}$ and $T_2 \subset \Z/{M_2}$
we recall that we have already defined their \lq intersection\rq\ 
$T_1 \cap T_2$ to be the set of 
all elements of $\Z/{{\rm lcm}(M_1,M_2)}$
whose reduction modulo $M_1$ and $M_2$ is respectively in $T_1$ and $T_2$.

The following Proposition 
will be our main tool in proving Proposition~\ref{prop:Fib733}.
\begin{prop}\label{prop:criterion}
Suppose $S=\left\{ (l_1,q_1), \ldots, (l_t,q_t) \right\}$ is a finite set of pairs
of primes $(l,q)$ where each $l$ satisfies the condition~(\ref{eqn:condl})
and each $q$ is $\geq 5$.
Let $K(S)={\rm lcm}(6,l_1 -1, \ldots, l_t -1)$, and  
\[
\N(S) = \cap_{(l,q) \in S} \N(l,q) \subset (\Z/K(S))^*.
\]
Write 
\[
\N(S) = \left\{ \overline{1}, \overline{a}, \overline{b},\ldots \right\} 
\quad \text{where} \quad 1< a<b< \cdots < K(S). 
\]
If $(n,y,p)$ is a solution to~(\ref{eqn:FibPrimeIndex}) with $p > q_1,\ldots,q_t$
then $n \geq a$.
\end{prop}
\begin{proof}
Suppose that $(n,y,p)$ is a solution to the equation~(\ref{eqn:FibPrimeIndex}). 
It follows immediately from Lemma~\ref{lem:Nl}  and the
definition of \lq intersection\rq\ that the reduction
of $n$ modulo $K(S)$ belongs to $\N(S)$. 

The reader can check for himself that $\overline{1}$
is always in $\N(S)$. Moreover $n \neq 1$ since 
$n > 25000$. Hence $n \geq a$.
\end{proof}
The following Lemma will provide a useful check for our later calculations.
\begin{lem}\label{lem:4elements}
With the notation of the above proposition, suppose that $4|K(S)$. Then
the residue classes of $1,-1,K(S)/2+1,K(S)/2-1$ modulo $K(S)$ all belong
to $\N(S)$.
\end{lem}
\begin{proof}
	We note that $H_1=H_{-1}=1$, and so $E_1=E$.
	It follows from the definition of $\N(l,p)$ that the residue classes of $1,~-1$
	modulo $K(l)$ belong to $\N(l,q)$ for all $(l,q) \in S$, and so  
	residue classes of $1,~-1$ modulo $K(S)$
	belong to $\N(S)$.

	Let us prove the same for $n=K(S)/2+1$; we will leave the other case to
	the reader. Suppose that $(l,q) \in S$.
	We would like to prove that the residue class of $n$ modulo
	$K(l)$ belongs to $\N(l,q)$. 
	Write $v_2 : \Z \rightarrow \Z_{\geq 0} \cup \left\{ \infty \right\}$
	for the $2$-adic valuation.

	Clearly $(l-1) \mid K(S)$. If $v_2(l-1) < v_2(K(S))$ then 
	$n \equiv 1 \pmod{K(l)}$ and we already know that 
	$\overline{1} \in \N(l,q)$. 

	Thus suppose that $v_2(l-1)=v_2(K(S))$. Since $4|K(S)$
	we see that $l \equiv 1 \pmod{4}$. Further we can write
	\[
		n=\frac{K(S)}{2}+1=k \frac{(l-1)}{2}+1
	\]
	for some odd integer $k$. Note that
	\[
		\omega^n \equiv (\omega^\frac{l-1}{2})^k \cdot \omega
		\equiv \pm \omega \pmod{l},
	\]
	and so
	\[
		H_n \equiv \pm L_n \equiv \pm (\omega^n -\omega^{-n})
		\equiv \pm (\omega - \omega^{-1}) \equiv  \pm 1 \pmod{l}.
	\]
	A glance at the definition of $\N(l,q)$
	shows that we must prove that $a_l(E_n) - a_l(E)$ is divisible by some prime
	greater than $q$. Actually we will 
	prove that $a_l(E_n)=a_l(E)$. By comparing the equations for $E$ and $E_n$
	we see that, modulo $l$, the curves $E$ and $E_n$ are isomorphic when $H(n) \equiv 1
	\pmod{l}$. If $H(n) \equiv -1 \pmod{l}$ then, modulo $l$, the curve $E_n$ is the
	$-1$-twist of $E$. But $\kro{-1}{l}=1$, and so again $E$ and $E_n$ are
	isomorphic modulo $l$. This proves  that $a_l(E_n)=a_l(E)$ and completes the
	proof.
\end{proof}
\subsection{Proof of Proposition~\ref{prop:Fib733}}
Suppose that $(n,y,p)$ is a solution to~(\ref{eqn:FibPrimeIndex}).
Notice that Proposition~\ref{prop:criterion} provides us with 
a way of obtaining lower bounds for $n$, and Proposition~\ref{prop:Bugeaud1}
provides us with a way of obtaining an upper bound for $n$ (dependent on $p$).
This gives us hope, given a particular prime $p$, that we may be
able to obtain a contradiction using these two Propositions and
so prove that there are no solutions for this particular $p$. 

We wrote a {\tt PARI/GP} program to carry out the above idea 
and derive the contradiction for the primes in the range $7 \leq p \leq 733$.

We would like to give the reader the flavour of this computation by providing
more for the proof that $p>7$.  

A priori, all we know about $p$ is that $p>5$, so we take $q=5$.
We let $S=\left\{ (11,5) \right\}$. Then
\[
\N(S) = \N(11,5) =\left\{ \overline{1}, \overline{11}, \overline{19}, \overline{29} \right\} \subset \Z/30,
\]
where we used our program to calculate $\N(11,5)$ from the definition of $\N(l,q)$.
Next we look for primes $l$ satisfying $l \equiv \pm 1 \pmod{5}$ and 
\[
(l-1) | M, \quad \text{where } M=6983776800= 2^5 \times 3^3 \times 5^2 \times 7 \times 11 \times \ldots \times 19  
\]
and for each such prime $l$ we find we append $(l,5)$ to the set $S$,
thus redefining $\N$ to be
$\N(S) \cap \N(l,5)$. We continue until $\N \subset \Z/M$
and $\N(S)$ has $4$ elements (we do not expect less than $4$ elements
by Lemma~\ref{lem:4elements}).
The reader will no doubt expect that since most of our $l-1$ are highly
composite and have lots of common factors, the set $\N(S)$ will be a small
set of congruences modulo a large modulus.
After a few seconds we found that
\[ 
\N(S) = \left\{ 
\overline{1}, \overline{3491888399}, \overline{3491888401}, 
\overline{6983776799}
\right\}
\subset \Z/M.
\]
We then replaced the value of $M$ by $M \times 23$ and continued until $\N(S)$ had
exactly $4$ elements and $\N \subset \Z/M$ with this new value of $M$, etc.
The entire computation took $42$ seconds and proved that the
the reduction of $n$ belongs to a set
\[
\N= \left\{
\overline{1},
\overline{a},
\overline{b},
\overline{c}
\right\} \subset \Z/M 
\]
where
\[
\begin{array}{l}
a=100704598854427777024179418273944411482999002799, \\
b=100704598854427777024179418273944411482999002801, \\ 
c=201409197708855554048358836547888822965998005599. \\
\end{array}
\] 
and the value of $M$ is now
\[
M=2^5 \times 3^3 \times 5^2 \times 7
\times 11 \times \ldots \times 109.  
\]
Note that $a \approx 1{.}007 \times 10^{47}$. 

By Proposition~\ref{prop:criterion}, $n \geq a$. However, if $p=7$ 
then Proposition~\ref{prop:Bugeaud1} implies that $n < 2{.}639 \times 10^{46}$. 
This proves that $p >7$.  As a check on our computations, we note that
$a=M/2-1$, $b=M/2+1$ and $c=M-1$ which is entirely consistent with 
Lemma~\ref{lem:4elements}.

The next step is to prove that $p >11$. We continue as above but
now take $q=7$. We note that $\N(l,7) \subseteq \N(l,5)$ for any
prime $l$, and that probably $\N(l,7)$ is strictly smaller $\N(l,5)$. Thus
our sieve becomes more efficient. 

The proof program took roughly 97 hours to run on a $1.7$ GHz Intel Pentium.
By the end of the proof the set $S$ had $6262$ pairs, and we have also
shown that $p > 733$ and $n \geq 1{.}033 \times 10^{8733}$. To complete the proof
we must show that $\log{y} \geq 10^{8000}$. However  
$y^p=F_n=(\omega^n-\tau^n)/\sqrt 5$. Taking logarithms and using
Peth\H{o}'s result  that $p < 5{.}1 \times 10^{17}$ (mentioned Section~\ref{sec:survey}) 
we deduce that $\log{y} \geq 10^{8000}$
with a huge margin.

\subsection{Proof of Proposition~\ref{prop:Luc283}}
The proof of Proposition~\ref{prop:Luc283} is
practically identical to the above proof of Proposition~\ref{prop:Fib733}
and we omit almost all the details. We can take $K(l)=l-1$ in this case,
and we let $E_n$ and $E$ be given by equations~(\ref{eqn:LucFrey})
and~(\ref{eqn:LucEll}) respectively. If $l \equiv \pm 1 \pmod{5}$ 
and $q \geq 5$ are primes we define  $\N(l,q)$ to be the subset of 
all $n \in (\Z/K(l))^*$ such that
\begin{itemize}
        \item either $5 F_n^2 -4 \not \equiv 0 \pmod{l}$,
        and the integer $a_l(E_n) - a_l(E)$ is divisible by some prime $p > q$.
        \item or $5 F_n^2 -4 \equiv 0 \pmod{l}$ and
        one of the two integers  $l+1 \pm a_l(E)$ is divisible by some prime
        $p > q$.
\end{itemize}
The other details are practically identical to the Fibonacci case. 
Since the lower bound for $p$
that we are trying to establish is much smaller in the Lucas case
our program runs much faster, and completes the proof on the same machine
in about $6$ hours.

 
\section{A Refined Bound on $p$ using Linear Form in $2$ Logarithms}
In the previous section we showed that if $(n,y,p)$ is solution to
equation~(\ref{eqn:LucPrimeIndex}) then $p > 283$.
In this section we will use the results of the paper of
Laurent, Mignotte and Nesterenko \cite{LMN} on linear forms
in $2$ logarithms to prove that $p \leq 283$ thus completing
the proof of Theorem~\ref{thm:Lucas}. 

The Fibonacci case still needs more work, since it yields a linear
form in 3 logarithms. However, for now we are able to show the
following.
\begin{prop}\label{prop:pbig}
If $(n,y,p)$ is a solution to equation~(\ref{eqn:FibPrimeIndex})
then $p > 2 \times 10^{8}$.
\end{prop}
\begin{proof}
Suppose that $(n,y,p)$ is a solution to~(\ref{eqn:FibPrimeIndex}).
%
The most obvious approach to obtain an upper bound for $p$ 
is to consider
\[
F_n=\frac{\omega^n-\omega^{-n}}{\sqrt{5}}=y^p
\]
and the linear form in logarithms
\[
\Lambda = n \log{\omega} - \log{\sqrt{5}} - p \log{y}.
\]
Then a standard argument shows that 
\[
\log |\Lambda| < -2p \log{y} +1.
\]
We note that $\Lambda$ is a linear form in $3$ logarithms.
In the remainder of this paper we will present a 
substantial improvement to the theory of linear form
in $3$ logarithms, and apply our result to show that $p < 2 \times 10^{8}$.

For now, to prove the proposition, we argue by contradiction, assuming that
$p < 2 \times 10^8$.  We then know from Proposition~\ref{prop:Fibnequivpm1}
that $n \equiv \pm 1 \pmod{p}$ for primes $p$ in this range.
Write $n=s p+\epsilon$, where $\epsilon=\pm 1$.
Note now that we can rewrite $\Lambda$ as
\[
\Lambda = p \log \,(\omega^s/y) - \log \, (\sqrt{5} \omega^{-\epsilon}),
\]
which is now a linear form in $2$ logarithms.  
We can apply Th\'eor\`eme~1 of \cite{LMN} with
\[
\Lambda = b_1 \log \alpha_1 - b_2 \log \alpha_2,
\]
where
\[
b_1=p, \quad \alpha_1 = \omega^s/y;
\qquad b_2 = 1, \quad  \alpha_2 = \sqrt 5 \,\omega^{-\epsilon}
\]
and
\[
\log \alpha_2  = \log \sqrt 5 - \epsilon \log\omega, \quad
\log \alpha_1 \approx \frac{1}{p}\,\log \alpha_2, 
\]
and
\[
\h (\alpha_2) = \frac{\log 10}{2}, \quad \h (\alpha_1) \le \log y + \log 5.
\]
Thus (with the notation of this result), we can take
\[
a_1 = (\rho-1) \log \alpha_1 + 4\,(\log y + \log 5), \quad
a_2 =  (\rho-1) \log \alpha_2 + 2\,\log 10.
\]
\bigskip

\noindent {\it The case $\epsilon=-1$}
\medskip

In this case, we choose (again with the notation of \cite{LMN})
$L=8$, $\rho=27{.}6$, $m=0{.}209671121$ and get
\[
p \le 733.
\]
\bigskip

\noindent {\it The case $\epsilon=1$}
\medskip

In this case, we choose 
$L=7$, $\rho=31{.}6$, $m=0{.}218149476$ and get
\[
p \le 241.
\]
In either case we have $p \leq 733$ which
contradicts Proposition~\ref{prop:Fib733}.
This completes the proof of the Proposition.
\end{proof}
As promised we also complete the resolution of the Lucas
case by presenting the proof of Theorem~\ref{thm:Lucas}.

\begin{proof}[Proof of Theorem~\ref{thm:Lucas}]

Suppose that $(n,y,p)$ is a solution to
equation~(\ref{eqn:LucPrimeIndex}).
It is apparent, by Proposition~\ref{prop:Luc283},
that all we have to do is to show that
$p \leq 283$, and to do this we apply
\cite{LMN}.

Put
\[
\Lambda = p \log y - n \log \omega
\]
where $\omega = (1+\sqrt 5)/2$. 
By Proposition~\ref{prop:Luc283} 
we know that 
\[
\log y >10^6,
\]
and indeed much more.
Then (because $L_n=\omega^n+(-1/\omega)^n$)
\[
\log | \Lambda | < - 2p\log y +1.
\]
Write $n=sp+r$ with $0\le r < p$. This allows us to
rewrite $\Lambda$ as
\[
\Lambda = p \log (y/\omega^s) - r \log \omega .
\]
We apply \cite[Proposition 1]{LMN} with the notation $D=2$ and
\[
\alpha_1=y/\omega^s, \
\alpha_2=\omega, \
b_1=p, \
b_2=r, \
a_1=2{.}00001\,\log (y\omega^s), \
a_2=(\rho+1)\log \omega.
\]
Here
\[
\tilde b = \frac{1}{\rho+1}\left( \frac{p}{\log \omega} + 
\frac{r(1+\rho)}{a_1}\right)
\approx \frac{1}{\rho+1}\, \frac{p}{\log \omega}.
\]
We get either
\[
p \le \mu L (\rho+1) \log \omega
\]
or
\[
\log | \Lambda | > - KL \log \rho -\log (KL)
\]
provided that
\[
2K\log \theta +2\log(2\pi K/e^{3/2})-3\log(KL)- \frac{c+\log K}{3K}
- \frac{a_1 L}{3} - \frac{a_2 L^2}{3}- \frac{2K}{\mu a_1}\ge 0
\]
and
\[
\mu \bigl( (L-1)\log \rho +2\log(2/\theta)-2(1{.}5-\log \mu +\log
\tilde b)\bigr)\ge
\frac{L}{3}.
\]
It is enough to take
\[
\log \theta = \frac{1{.}01 \times a_1 L}{3 \times 2 \mu^2 a_1 a_2 L}
= \frac{1{.}01}{6 \mu^2 (\rho +1) \log \omega}.
\]
For $\rho=22.9$, taking $\mu=2/(3\omega)$,
working as above we first get $p<326$ and then,
after several iterations of the above argument,
\[
p\le 283.
\]
\end{proof}


\section{An Estimate on Linear Forms in Three Logarithms}\label{sec:3logs} 

\subsection{Preliminaries}
\begin{lem}\label{lem:Ma1}
Let $K$, $L$, $R$, $S$, $T$ be positive integers, put $N=K^2L$ and
assume $N\le RST$, put also
\[
\ell_n =\left\lfloor \frac{n-1}{K^2}\right\rfloor, \quad 1\le n \le N,
\]
and $(r_1,\ldots,r_N)\in \{0,1,\ldots,R-1\}^N$. Suppose that for each
$r\in
\{0,1,\ldots,R-1\}$ there are at most $ST$ indices such that $r_j=r$.
Then
\[
\left| \sum_{n=1}^N \ell_n r_n -M\right|\le G_R
\]
where
\[
M=\left( \frac{L-1}{2}\right) \sum_{n=1}^N r_n
\qquad {\rm and}\qquad
G_R=\frac{NLR}{2} \left(\frac{1}{4}- \frac{N}{12RST}\right).
\]
\end{lem}
\begin{proof}
Apply \cite[Lemme 4]{LMN}. 
\end{proof}

As in \cite{B} or \cite[page 192]{W}, for $(k,m)\in \bN^2$, we put 
$\Vert (k,m)\Vert=k+m$. And we put
\[
\Theta (K_0,I)=\min \bigl\{ \Vert (k_1,m_1)\Vert+ \cdots +\Vert
(k_I,m_I)\Vert \bigr\},
\]
where the minimum is taken over if the $I$--tuples $(k_1,m_1)$, \dots,
$(k_I,m_I)\in \bN^2$ which are pairwise distinct and satisfy $m_1$,
\dots, $m_I\le K_0$.
Then, we have:

\begin{lem}\label{lem:Ma2}
Let $K_0$, $L$ and $I$ be positive integers with $K_0\ge 3$, $L\ge 2$
and
$I\ge K_0(K_0+1)/2$. Then
\[
\Theta (K_0,I)\ge
\left(\frac{I^2}{2(K_0+1)}\right)
\left(1+  \frac{(K_0-1)( K_0+1)}{I} - \frac{K_0(K_0+2)(K_0+1)^2}{12 I^2} 
\right).
\]
\end{lem}
\begin{proof}
Except for some details, this is \cite[Lemma 1.4]{B}. We follow more or
less the proof of this result.
The argument is elementary: the smallest value for the sum
$\Vert (k_1,m_1)\Vert + \cdots + \Vert (k_I,m_I)\Vert$ is reached when
we choose
successively, for each integer $n=0$, 1, \dots\ all the points in the
domain
\[
D_n = \{(k,m)\in \bN^2;\; m\le K_0, \; k+m=n\},
\]
and stop when the total number of points is $I$. Moreover,
\[
{\rm Card} (D_n) =
\left\{
\begin{array}{ll}
n+1, & \text{if $n\le K_0$}, \\
K_0+1, & \text{if  $n\ge K_0$.}
\end{array}
\right.
\]
With the notation of \cite{B}, the number $I$ of points can be written as
\[
I=\left(A-\frac{K_0}{2}\right) (K_0+1)+r, \quad {\rm with}\ \; 0\le r
\le K_0,
\]
provided that $I\ge K_0(K_0+1)/2$, which is a hypothesis of the
Lemma.
\smallskip
Then, the computation of \cite{B} shows that
\[
\Theta (K_0,I)
\ge \tilde\Theta (K_0,I):=
\frac{K_0+1}{2}\, 
\left(A(A-1)
-\frac{K_0(K_0-1)}{3} \right)+rA.
\]
In terms of $I$,
\[
A=\frac{K_0}{2}+\frac{I-r}{K_0+1}.
\]
We have,
\[
\frac{\partial \tilde \Theta}{\partial r }
= \frac{K_0+1}{2}\,  (2A-1) \frac{\partial A}{\partial r }+ A
+ r \frac{\partial A}{\partial r } = - \frac{2A-1}{2}+A- \frac{r}{K_0+1}
=\frac{1}{2}- \frac{r}{K_0+1},
\]
which shows that the minimum of $\tilde \Theta$ is reached either for
$r=0$ or $r=K_0$. It is
easy to verify that $\tilde \Theta$ takes the same value for $r=0$ and
$r=K_0+1$ (which is
indeed out of the range of $r$), this implies that the minimum is
reached for $r=0$. It
follows that
\[
\begin{split}
\frac{2\Theta (K_0,I)}{K_0+1}
& \ge
\left(\frac{K_0}{2}+\frac{I}{K_0+1} \right)\left(\frac{K_0}{2}+
\frac{I}{K_0+1}-1 \right)
- \frac{K_0(K_0-1)}{3}
\\
& = \frac{K_0^2}{4}+\frac{I^2}{(K_0+1)^2} + \frac{K_0 I}{K_0+1}
- \frac{K_0}{2}- \frac{I}{K_0+1}
- \frac{K_0^2}{3}+ \frac{K_0}{3}
\\
&= \frac{I^2}{(K_0+1)^2}+\frac{(K_0 -1)I}{K_0+1}-
\frac{K_0^2}{12}- \frac{K_0}{6}
\\
&=
\left(\frac{I}{K_0+1}\right)^2
\left(1+  \frac{(K_0-1)( K_0+1)}{I} - \frac{K_0(K_0+2)(K_0+1)^2}{12 I^2} 
\right),
\end{split}
\]
which proves the  lemma. 
\end{proof}

The version of Liouville inequality that we use is the same as in \cite[pages 298--99]{LMN}:

\begin{lem}\label{lem:Ma3}
Let $\alpha_1,~\alpha_2,~\alpha_3$ be non-zero algebraic numbers
and
$f \in \Z[X_1,X_2,X_3]$
such that $f(\alpha_1,\alpha_2,\alpha_3)\not=0$, then
\[
\begin{split}
| f(\alpha_1,\alpha_2,\alpha_3)| \ge 
& |f|^{-{D+1}} (\alpha_1^*)^{d_1}
(\alpha_2^*)^{d_2}(\alpha_3^*)^{d_3} \\
& \times \exp \bigl\{-D\bigl( d_1 \h(\alpha_1)+d_2 \h(\alpha_2)+d_3
\h(\alpha_3)\bigr) \bigr\}
\end{split}
\]
where
$\,D=[\Q(\alpha_1,\alpha_2,\alpha_3) :
\Q]\bigm/[\R(\alpha_1,\alpha_2,\alpha_3) : \R]$,
\[
d_i=\deg_{X_i}f, \ \  i=1,2,3,
\qquad |f|=\max \bigl\{|f(z_1,z_2,z_3)|;\; |z_i|\le 1, \
i=1,2,3\bigr\},
\]
and $\h(\alpha)$ is the absolute logarithmic height of the algebraic
number $\alpha$, and
$\alpha^*=\max\{1,|\alpha|\}$.
\end{lem}
\begin{lem}\label{lem:Ma4}
let $K> 1$ be an integer, then
\[
\log \, \left(\prod_{k=1}^{K-1}k!\right)^\frac{4}{K(K-1)}
\ge 2\log K - 3+ \frac{2\log (2\pi K/e^{3/2})}{K-1}-
\frac{2+6\pi^{-2}+\log K}{3K(K-1)}.
\]
\end{lem}
\begin{proof}
This is a consequence of a variant of the proof of \cite[Lemme 8]{LMN}.
\end{proof}

Now we present the type of linear forms in three logs that we shall
study. For a while, we consider
three non-zero  algebraic numbers $\alpha_1$, $\alpha_2$ and $\alpha_3$
and positive
rational integers $b_1$, $b_2$, $b_3$ with $\gcd(b_1,b_2,b_3)=1$, and
the linear form
\[
\Lambda = b_2 \log \alpha_2-b_1\log \alpha_1-b_3\log \alpha_3\not=0,
\]
without any loss in generality.

We restrict our study to the following cases,
\begin{itemize}
\item the real case: $\alpha_1$, $\alpha_2$ and $\alpha_3$ are real numbers
${}>1$, and the logarithms
of the $\alpha_i$'s are real (and ${}>0$),
\item the complex case: $\alpha_1$, $\alpha_2$ and $\alpha_3$ are complex
numbers of modulus one, and the logarithms
of the $\alpha_i$'s are arbitrary determinations of the logarithm.
\end{itemize}
This does not cause inconvenience in practice since in the general case
we obviously always have 
\[
|\Lambda| \ge \max\{|\Re(\Lambda)|, |\Im(\Lambda)|\}.
\]

Without loss of generality, we may assume that
\[
b_2 |\log \alpha_2|=  b_1|\log \alpha_1|+b_3|\log \alpha_3|\pm
|\Lambda|.
\]
We shall choose rational positive integers $K$, $L$, $R$, $S$, $T$,
with $K$, $L\ge 2$,
we put $N=K^2L$ and we assume $RST\ge N$.
Let $a_1$, $a_2$, $a_3$ be positive real numbers, without loss of
generality, let $a_1\ge a_3$.

The authors of \cite{B} use Laurent's method, and they consider a suitable
interpolation determinant $\Delta$.
Let $i$ be an index such that $(k_i,m_i,\ell_i)$ runs trough all
triples of integers with
$0\le k_i\le K-1$, $0\le m_i\le K-1$ and $0\le \ell_i\le L-1$. So each
number $0$, \dots, $K-1$
occurs $KL$ times as a $k_i$, and similarly as a $m_i$, and each number
$0$, \dots, $L-1$
occurs $K^2$ times as an $\ell_i$.
With the above definitions, let
\[
\Delta = \det
\left\{ \binom{r_jb_2+s_jb_1}{k_i} \binom{t_jb_2+s_jb_3}{m_i}
\alpha_1^{\ell_i r_j}\alpha_2^{\ell_i s_j}\alpha_3^{\ell_i t_j}
\right\}
\]
where $r_j$, $s_j$, $t_j$ are non-negative integers less than $R$, $S$,
$T$, respectively,
such that $(r_j,s_j,t_j)$ runs over $N$ distinct triples.
Put $\,\beta_1=b_1/b_2$, $\beta_3=b_3/b_2$. Let
\[
\lambda_i = \ell_i - \frac{L-1}{2},\quad
\eta_0 = \frac{R-1}{2}+\beta_1 \frac{S-1}{2},\quad
\zeta_0 = \frac{T-1}{2}+\beta_3 \frac{S-1}{2},\quad
\]
and
\[
b = (b_2\eta_0) (b_2\zeta_0)\left(\prod_{k=1}^{K-1}k!\right)^{-\frac{4}{K(K-1)}}.
\]
Notice that, by Lemma~\ref{lem:Ma4},
\[
\begin{split}
\log b
&\le \log \frac{(R-1)b_2+(S-1)b_1}{2}+\log \frac{(T-1)b_2+(S-1)b_3}{2}
\\
&\quad -2\log K + 3- \frac{2\log (2\pi K/e^{3/2})}{K-1}+
\frac{2+6\pi^{-2}+\log K}{3K(K-1)}.
\end{split}
\]

Then $\,\sum_{i=0}^{N-1}\lambda_i=0\,$ and (\cite[formula (2.1)]{B})
\[
\alpha_1^{\ell_i r_j}\alpha_2^{\ell_i s_j}\alpha_3^{\ell_i t_j}
= \alpha_1^{\lambda_i(r_j+s_j\beta_1)}
\alpha_3^{\lambda_i(t_j+s_j\beta_3)}(1+\theta_{ij}\Lambda'),
\]
where
\[
\Lambda'=|\Lambda| \cdot
\max \left\{
\frac{LRe^{LR |\Lambda|/(2b_1)}}{2b_1},
\frac{LSe^{LS |\Lambda|/(2b_2)}}{2b_2},
\frac{LTe^{LT |\Lambda|/(2b_3)}}{2b_3}
\right\}
\]
and where all $|\theta_{ij}|$ are ${}\le 1$.

\subsection{An upper bound for $|\Delta|$}
It is proved in \cite{B} (last formula of page~111) that
\[
\Delta = \alpha_1^{ M_1} \alpha_2^{ M_2}\alpha_3^{ M_3}
\sum_{\cI \subseteq \cN}(\Lambda')^{N-|\cI|}\Delta_\cI
\]
where
\[
M_1=\frac{L-1}{2}\sum_{j=1}^N r_j, \qquad
M_2=\frac{L-1}{2}\sum_{j=1}^N s_j, \qquad
M_3=\frac{L-1}{2}\sum_{j=1}^N t_j, 
\]
and where
$\cN=\{0,1,\ldots,N-1\}$ and $\Delta_\cI$ is the determinant of a
certain matrix $\cM_\cI$
defined below. Let
\[
\phi_j(z,\zeta)=\frac{b_2^{k_i+m_i}}{k_i!\,m_i} z^{k_i}\zeta^{m_i}
\alpha_1^{\lambda_i z}\alpha_3^{\lambda_i \zeta},
\]
[where $\,\alpha_1^{\lambda_i z}=\exp(\lambda_i z \log \alpha_1)\,$
and similarly for $\,\alpha_3^{\lambda_i \zeta}$]
and
\[
\Phi_\cI (x)_{ij}=
\begin{cases}
\phi_j(xz_j,x\zeta_j) & \text{if $\,i\in \cI$,} \\
\theta_{ij}\phi_j(xz_j,x\zeta_j) & \text{if $\,i\not\in \cI$.}
\end{cases}
\]
Then, $\cM_\cI = \bigl(\Phi_\cI(1)_{ij}\bigr)$ and letting 
$\Psi_\cI (x)=\det \bigl( \Phi_\cI (x)\bigr)$, gives
\[
|\Delta_\cI | = \left|\det \bigl( \Phi_\cI (1)\bigr)\right| = |\Psi_\cI(1)|.
\]
Now, let
\[
J_\cI = {\rm order} (\Psi, 0),
\]
the maximum modulus principle implies
\[
|\Psi_\cI(1)|\le \rho^{-J_\cI}\cdot \max_{|x|=\rho }|\Psi_\cI(x)|.
\]
Since $|z_j|\le \eta_0$ and $|\zeta_j|\le \zeta_0$,
\[
\begin{split}
\max_{|x|=\rho} \bigl|\Psi_\cI (x) \bigr| \le & 
N!\,
\frac{{b_2}^{\sum k_i+\sum m_i}}{\prod k_i!\, \prod m_i!}  \;
(\rho \eta_0)^{ \sum k_i } (\rho \zeta_0)^{ \sum m_i }\!
\\
& \times \max_{\sigma\in \gS (\cN)} 
  \exp
\left\{
 \rho \Bigl(
\bigl (\sum \lambda_i z_{\sigma(i)} \bigr) \log \alpha_1 +
 \bigl(\sum \lambda_i \zeta_{\sigma(i)} \bigr) \log \alpha_2
\Bigr)
\right\}.
\end{split}
\]
Put
\[
g=\frac{1}{4}- \frac{N}{12RST}, \quad
G_1=\frac{NLR}{2} \,g, \quad
G_2=\frac{NLS}{2} \,g, \quad
G_3=\frac{NLT}{2} \,g,
\]
then (see the proof of \cite{B} p.~114 and use Lemma~\ref{lem:Ma1})
\[
\sum_{i=0}^{N-1} \lambda_i z_{\sigma(i)} \le G_1+\beta_1 G_2,
\qquad
\sum_{i=0}^{N-1} \lambda_i \zeta_{\sigma(i)} \le G_3+\beta_3 G_2.
\]
It follows that
(recall that $b_2 |\log \alpha_2|=  b_1|\log \alpha_1|+b_3|\log
\alpha_3|\pm |\Lambda|$)
\[
\begin{split}
& \exp
\left\{
\rho
 \Bigl(
   \bigl (\sum \lambda_i z_{\sigma(i)} \bigr)|\log \alpha_1 |+
   \bigl(\sum \lambda_i \zeta_{\sigma(i)} \bigr)|\log \alpha_3|
\Bigr)
\right\}
\\
 \le &
\exp
\left\{
\rho
  \Bigl(
 (G_1+\beta_1 G_2)|\log \alpha_1| +
 (G_3+\beta_3 G_2)|\log \alpha_3|
\Bigr)
\right\}
\\
 \kern-6mm\le &
\exp
\left\{
 \rho \left(
  G_1 |\log \alpha_1 |+ G_2 \Bigl(|\log \alpha_2|+\frac{|\Lambda|}{b_2}\Bigr)+
  G_3 |\log \alpha_3|
\right)
\right\}.
\end{split}
\]
As in \cite{B}, we see that if
\begin{equation}\label{eqn:*}
\Lambda' < \rho^{-KL}
\end{equation}
then
\[
\rho G_2 \, \frac{|\Lambda|}{b_2}
\le \frac{\rho K^2 L}{4\rho^{KL}}
\le \frac{e K^2 L}{4 e^{KL}}
\le \frac{ K^2 L^2}{e^{KL}} <10^{-4}
\]
for $KL\ge 15$. Putting these estimates together, we get that condition~(\ref{eqn:*}) implies
the upper bound
\[
\begin{split}
|\Delta| \le 1{.}0001 \,  {\alpha_1}^{M_1+\rho G_1}
\, {\alpha_2}^{M_2+\rho G_2}\, {\alpha_3}^{M_3+\rho G_3}
\,N! &\times 2^N \, \rho^{\,\sum (k_i+m_i)}
\\
\times \frac{(b_2\eta_0)^{\,\sum k_i }}{\prod k_i!}
\times \frac{(b_2\zeta_0)^{\,\sum m_i }}{\prod m_i!}
& \times \max_{\sigma\in \gS (\cN)} \, \frac{|\Lambda'|^{N-|\cI|}}{\rho^{J_\cI}}
\end{split}
\]
where
\[
J_\cI ={\rm order} (\Psi_\cI,0).
\]
Under condition~(\ref{eqn:*}), we have
\[
\frac{|\Lambda'|^{N-|\cI|}}{\rho^{J_\cI}}
\le \rho^{-KL(N-|\cI|)-J_\cI}.
\]
If $\,|\cI|\le 0{.5}\,N\,$ then
\[
KL(N-|\cI|)\ge 0{.5}\,KLN \ge \frac{NKL}{4}\left(1+\frac{4}{L}+\frac{1}{2K-1}\right)
\]
as soon as $\,K\ge 3\,$ and $\,L\ge 5$, conditions that we assume from
now on.

If $\,|\cI|\ge 0{.5}\,N$, then using \cite[Lemma 1.3]{B}, we obtain
\[
J_\cI \ge \Theta (K_0,|\cI|),
\qquad \text{for \ $\,K_0=2(K-1)\,$.}
\]
Now, $\,|\cI|\ge 0{.5}\,K^2L\,$ implies $\,|\cI|\ge 2{.5}\,K^2\,$
 and using  Lemma~\ref{lem:Ma2} we get (with the notation $I=|\cI|$)
\begin{multline*}
KL(N-I)+J_\cI  \ge 
KL(N-I) \\
 + \left( \frac{I^2}{2(K_0+1)}\right)
\left(1+  \frac{(K_0-1)( K_0+1)}{I} - \frac{K_0(K_0+2)(K_0+1)^2}{12 I^2} 
\right).
\end{multline*}
It is easy to verify that the right-hand side is a decreasing function
of $I$ in the range
$[N/2,N]$, since $L\ge 5$, and we get (recall that $N=K^2L$ and
$K_0=2K-2$)
\[
\begin{split}
 KL(N-|\cI|)+J_\cI
 &\ge
  \frac{N^2}{2(K_0+1)}
\left(1+  \frac{K_0^2-1}{N} - \frac{K_0 (K_0+2)(K_0+1)^2}{12 N^2} 
\right)
\\
&=  \frac{N^2}{4K}
\left( \frac{2K}{K_0+1}+ \frac{2K(K_0 -1)}{N} - \frac{K K_0(K_0+1)(K_0+2)}{6 N^2}  \right)
\\
&=  \frac{N^2}{4K}
\left(1+ \frac{1}{2K-1} + \frac{2(2K-3)}{KL} - \frac{ 2(K-1)(2K-1)}{3 K^2L^2}  \right)
\\
&=  \frac{N^2}{4K}
\left(1+ \frac{4}{L} + \frac{1}{2K-1}- \frac{4}{3L^2}
- \frac{6}{KL}+  \frac{2}{KL^2}- \frac{2}{3K^2L^2}   \right)
\\
&\ge  \frac{N^2}{4K}
\left(1+ \frac{4}{L} + \frac{1}{2K-1}- \frac{4}{3L^2} - \frac{6}{KL}   
\right),
\end{split}
\]
because $L\ge 5$,
and this implies, in all cases,
\[
KL(N-|\cI|)+J_\cI \ge
  \frac{N^2}{4K}
\left(1+  \frac{4}{L} + \frac{1}{2K-1}- \frac{6}{KL} - \frac{4}{3L^2}
\right).
\]
Thus, gathering all the previous estimates and using the relations
\[
\sum_{i=0}^{N-1}k_i = 
\sum_{i=0}^{N-1}m_i = 
\frac{(K-1)K}{2}\, KL=\frac{N}{2}\,(K-1),
\]
and the definition of $b$,
we obtain the following result.
\begin{prop}\label{prop:Ma1}
With, the previous notation, if  $K\ge 3$, $L\ge 5$ and $\Lambda'\le
\rho^{-KL}$,
with $\rho>1$,
then
\[
\begin{split}
\log |\Delta|
&\le  \sum_{i=1}^3  M_i \log |\alpha_i|+\rho \sum_{i=1}^3 G_i |\log
\alpha_i|
 +\log(N!)+N\log 2+ \frac{N}{2}\,(K-1) \log b
\\
& \kern 15mm 
-\left(
\frac{NKL}{4}+ \frac{NKL}{4(2K-1)} -\frac{NK}{3L}- \frac{N}{2}
\right)
\log \rho + 0{.}0001 .
\end{split}
\]
\end{prop}

\subsection{A lower bound for $|\Delta|$}

Using our Liouville estimate (Lemma~\ref{lem:Ma3}) and arguing as in \cite{B}, or
\cite[Lemme 6]{LMN}, we get the following.
\begin{prop}\label{prop:Ma2}
If $\,\Delta \not= 0\,$ then
\[
\begin{split}
\log |\Delta|
\ge &
- \frac{D-1}{2}\,N\log N \\
& +\sum_{i=1}^3 ( M_i +G_i)\log |\alpha_i|
-2D \sum_{i=1}^3  G_i \h(\alpha_i) - \frac{D-1}{2}(K-1)\,N\log b.
\end{split}
\]
\end{prop}
\begin{proof}
We have $\Delta=P(\alpha_1,\alpha_2,\alpha_3)$ where $P\in
\Z[X_1,X_2,X_3]$ is given by
\begin{multline*}
P(X_1,X_2,X_3) = \\
\sum_{\sigma\in \gS_N} {\rm sg}(\sigma) \cdot
\prod _{i=1}^N \binom{r_{\sigma(i)} b_2+s_ {\sigma(i)} b_1}{k_i}
\binom{t_{\sigma(i)} b_2+s_ {\sigma(i)} b_3}{m_i}
X_1^{n_{r\sigma}}X_2^{n_{s\sigma}}X_3^{n_{t\sigma}},
\end{multline*}
and where
\[
n_{r\sigma}=\sum_{i=1}^N \ell_i r_{\sigma(i)}, \quad
n_{s\sigma}=\sum_{i=1}^N \ell_i s_{\sigma(i)}, \quad
n_{t\sigma}=\sum_{i=1}^N \ell_i t_{\sigma(i)}.
\]
By Lemma~\ref{lem:Ma1},
\[
\bigl|\deg_{X_i} P -M_i\bigr|\le G_i, \quad \text{$i=1$, 2, 3.}
\]
Let
\[
V_i = \lfloor M_i+G_i\rfloor, \qquad U_i = \lceil M_i-G_i\rceil,
\quad \text{$i=1$, 2, 3,}
\]
then
\[
\Delta = {\alpha_1}^{V_1} {\alpha_2}^{V_2} {\alpha_3}^{V_3} \,
\tilde P(\alpha_1^{-1},\alpha_2^{-1},\alpha_3^{-1}),
\]
where
\[
 \deg_{X_i} \tilde P  \le V_i-U_i, \quad \text{$i=1$, 2, 3.}
\]
By our Liouville estimate
\[
\log \bigl| \tilde P(\alpha_1^{-1},\alpha_2^{-1},\alpha_3^{-1}) \bigr|
\ge - (D-1)\log |\tilde P|-D\sum_{i=1}^3 (V_i-U_i)\h(\alpha_i).
\]
Now we have to find an upper bound for $|\tilde P|$
(or for $|P|$ which is equal to $|\tilde P|$). By the multilinearity of
the determinant, for all
$\eta$, $\zeta\in \bC$,
\[
P(z_1,z_2,z_3)=
\det\left(
\frac{( r_j b_2+s_jb_1-\eta)^{k_i}}{k_i!}\frac{( t_j b_2+s_jb_3-\zeta)^{m_i}}{m_i!}
\cdot {z_1}^{\ell_i r_j} \cdot {z_2}^{\ell_i s_j}\cdot {z_3}^{\ell_i
t_j}
\right).
\]
Choose
\[
\eta = \frac{ (R-1)b_2+(S-1)b_1}{2}, \quad
\zeta = \frac{ (T-1)b_2+(S-1)b_3}{2}
\]
and notice that, for $1\le j\le N$,
\[
| r_j b_2+s_jb_1-\eta|^{k_i} \le \left(\frac{ (R-1)b_2+(S-1)b_1}{2}\right)^{k_i},
\]
\[
| t_j b_2+s_jb_3-\zeta|^{k_i} \le \left(\frac{ (T-1)b_2+(S-1)b_3}{2}\right)^{m_i}
\]
and that
\[
\sum_{i=0}^{N-1}k_i = \sum_{i=0}^{N-1}m_i = \frac{(K-1)K}{2}\,KL=\frac{N}{2}\,(K-1),
\]
then Hadamard's inequality implies
\begin{multline*}
|P| \le 
N^{N/2} \left(\frac{ (R-1)b_2+(S-1)b_1}{2}\right)^\frac{(K-1)N}{2}
\left(\frac{ (T-1)b_2+(S-1)b_3}{2}\right)^\frac{(K-1)N}{2} \\
\times
 \left(\prod_{i=0}^{K-1} k_i!\right)^{-1}  \left(\prod_{i=0}^{K-1}
m_i!\right)^{-1}.
\end{multline*}
Recall that
\[
b = (b_2\eta_0) (b_2\zeta_0)\left(\prod_{k=1}^{K-1}k!\right)^{-\frac{4}{K(K-1) }},
\]
where
\[
\eta_0 = \frac{R-1}{2}+\beta_1 \frac{S-1}{2},\qquad
\zeta_0 = \frac{T-1}{2}+\beta_3 \frac{S-1}{2}.
\]
Thus we get,
\[
|P| \le N^{N/2}\,b^{(K-1)N/2}.
\]
Collecting all the above estimates, we find
\[
\begin{split}
\log |\Delta| \ge & 
- (D-1)\,\left(
 \log \left(N^{N/2}\right) 
 + \frac{(K-1)N}{2}\,\log b
\right) \\
& - D \sum_{i=1}^3 (V_i-U_i) \h(\alpha_i)
 +\sum_{i=1}^3 V_i\log |\alpha_i|.
\end{split}
\]
The inequalities $\, D\h(\alpha_i)\ge \log |\alpha_i|\ge 0\,$ imply 
\[
V_i\log |\alpha_i|- D   (V_i-U_i) \h(\alpha_i)
\ge (M_i+G_i)\log |\alpha_i|-2D  G_i \h(\alpha_i)
\]
and the result follows.
\end{proof}

\subsection{Synthesis}
Under the hypotheses of Proposition~\ref{prop:Ma1} and Proposition~\ref{prop:Ma2}, we get
\begin{multline*}
- \frac{D-1}{2}\,N\log N
 +\sum_{i=1}^3 ( M_i +G_i)\log |\alpha_i|
-2D \sum_{i=1}^3  G_i \h(\alpha_i)
- \frac{D-1}{2}(K-1)\,N\log b
\\
\le  \sum_{i=1}^3  M_i \log |\alpha_i|+\rho \sum_{i=1}^3 G_i |\log
\alpha_i|
 +\log(N!)+N\log 2+ \frac{N}{2}\,(K-1) \log b
\\
 \kern 13mm
-\left(
\frac{NKL}{4}+\frac{NKL}{4(2K-1)} - \frac{NK}{3L}- \frac{N}{2}
\right)
\log \rho + 0{.}0001
\end{multline*}
Or, after some simplification,
\begin{multline*}
- \frac{D-1}{2}\,N\log N
 \le
\sum_{i=1}^3 G_i \bigl(\rho|\log \alpha_i|- \log |\alpha_i|+2D
\h(\alpha_i)\bigr) 
 +\log(N!) \\
+N\log 2 + \frac{K-1}{2} \,DN\log b
  -\left(  \frac{NKL}{4}+ \frac{NKL}{4(2K-1)}- \frac{KN}{3L}-
\frac{N}{2}\right)
\log \rho   + 0{.}0001 .
\end{multline*}
This result implies (divide by $N/2$ and use $N! < N(N/e)^N$, true for
$N>7$).
\begin{prop}\label{prop:Ma3}
With, the previous notation, if  $K\ge 3$, $L\ge 5$,
 $\rho>1$, and if $\Delta\not=0$ then
\[
\Lambda' > \rho^{-KL}
\]
provided that
\[
\begin{split}
\left(\frac{KL}{2}+ \frac{L}{4}-1  - \frac{2K}{3L} \right)  \log \rho 
 \ge
 & \, (D+1)\log N + gL(a_1R+a_2S+a_3T) \\
 & +D (K-1) \,  \log b - 2 \log(e/2),
\end{split}
\]
where the $a_i$'s satisfy
\[
a_i\ge\rho|\log \alpha_i|- \log |\alpha_i|+2D \h(\alpha_i), \qquad
\text{$i=1$, 2, 3.}
\]
\end{prop}

\begin{remark}
We notice that the statement of Proposition~\ref{prop:Ma3} is perfectly symmetric
with respect to the
$b_i$'s or the $\alpha_i$'s, except for the choice of $b$. From now on
we do not assume
that $b_1$ and $b_3$ are positive, but we still suppose that $b_2>0$
and that
\[
b_2 |\log \alpha_2| =  |b_1\,\log \alpha_1|+ |b_3 \,\log \alpha_3|\pm
|\Lambda|.
\]
\end{remark}
\subsection{Row rank}
To conclude we need to find conditions under which one of our
determinants $\Delta$ is non-zero,
a so called {\it zero lemma}.
We quote \cite[Theorem~3]{B} with some minor technical changes.
\begin{prop}\label{prop:Ma4}
Let $K$, $L$, $R$, $R_1$, $R_2$, $S$, $S_1$, $S_2$, $T$, $T_1$, $T_2$
be rational integers
all ${}\ge3$, with $K\ge2L$, $R>R_1+R_2$, $S>S_1+S_2$, $T>T_1+T_2$ and
$T_1\ge R_1$.
Let $b_1$, $b_2$, $b_3$ and $\alpha_1$, $\alpha_2$, $\alpha_3$ as above
and moreover assume that
$\alpha_1$, $\alpha_2$, $\alpha_3$ are multiplicatively independent.
If
\[
4(R_1+1)(S_1+1) \ge T_1+1,
\leqno{\rm (i)}
\]
\[
4(R_1+1)(T_1+1) \ge S_1+1,
\leqno{\rm (i)}
\]
\[
 (R_2+1)(S_2+1)(T_2+1) \ge 12(K-1)^2(L-1),
\leqno{\rm (iii)}
\]
and
\[
(R_1+1)(S_1+1)(T_1+1) \ge 8(2K+L-2)^2
\leqno{\rm (iv)}
\]
then {\bf either} there exists a choice of $\Delta$ which is non-zero
{\bf or}
 at least one of the following conditions hold
\[
 \exists r,\, s \in \Z, \quad \text{$rb_2=sb_1\;$ with $\;0<r\le
R_i\;$ and
 $\;0<s\le S_i\;$ for some $\;i=1$, 2,}
 \leqno{\bf (C1)}
\]
\[
 \exists t,\, s \in \Z, \quad \hbox{$tb_2=sb_3\;$ with $\;0<t\le
T_i\;$ and
 $\;0<s\le S_i\;$ for some $\;i=1$, 2,}
 \leqno{\bf (C2)}
\]
\noindent {\bf (C3)}: there exists $r'$, $s'$, $t'$, $t'' \in \Z$ such that
\[
 s't'b_1+r't''b_2+r's'b_3=0
\]
 which satisfy
\[
 0<|r'| <\min\left\{R_1+1,
 \left( \frac{(R_1+1)(S_1+1)}{T_1+1 } \right)^{1/2}\right\},
\]
\[ 
 0<|s'|<\min\left\{S_1+1,
 \left( \frac{(R_1+1)(S_1+1)}{T_1+1 } \right)^{1/2}\right\}
\]
 and
\[
  0<|t'|<\min\left\{T_1+1,
 \left( \frac{(S_1+1)(T_1+1)}{R_1+1 } \right)^{1/2}\right\},
\]
\[
\kern4mm |t''|<\min\left\{T_1+1,
 \left( \frac{(R_1+1)(T_1+1)}{S_1+1} \right)^{1/2}\right\}, 
\]
 which implies a non-trivial relation of the form
\[
 d_1b_1+d_2b_2+d_3b_3=0
 \qquad
 \text{with \ $|d_1|\le S_1$, $\,|d_2|\le R_1$,
 $\,|d_3|\le \frac{(R_1+1)(S_1+1)}{T_1+1}$.}
\]
\end{prop}
 
\subsection{A lower bound for the linear form}
Now we have all the tools to conclude and we get at once the following
result.

\begin{theorem}\label{thm:Maurice}
We consider three non-zero  algebraic numbers $\alpha_1$, $\alpha_2$
and $\alpha_3$
which are multiplicatively independent and positive
rational integers $b_1$, $b_2$, $b_3$ with $\gcd(b_1,b_2,b_3)=1$, and
the linear form
\[
\Lambda = b_2 \log \alpha_2-b_1\log \alpha_1-b_3\log \alpha_3\not=0.
\]
Where {\bf either}  $\alpha_1$, $\alpha_2$ and $\alpha_3$ are real
numbers ${}>1$, and the logarithms
of the $\alpha_i$'s are real (and ${}>0$),
{\bf or} $\alpha_1$, $\alpha_2$ and $\alpha_3$ are complex numbers of
modulus one, and the logarithms
of the $\alpha_i$'s are arbitrary determinations of the logarithm.
Without loss of generality, we assume that
\[
b_2 |\log \alpha_2|=  b_1|\log \alpha_1|+b_3|\log \alpha_3|\pm
|\Lambda|.
\]
 Let $K$, $L$, $R$, $R_1$, $R_2$, $S$, $S_1$, $S_2$, $T$, $T_1$, $T_2$
be rational integers
all ${}\ge3$, with
 $K\ge2L$, $L\ge 5$, $R>R_1+R_2$, $S>S_1+S_2$, $T>T_1+T_2$ and $T_1\ge
R_1$.
Let $\rho>1$ be a real number. Assume first that

\begin{eqnarray*}
{\rm (o)} \qquad
\left( \frac{KL}{2}+ \frac{L}{4}-1 - \frac{2K}{3L} \right)  \log \rho  
 \ge &
(D+1)\log N  + gL(a_1R+a_2S+a_3T) \\
 &+D (K-1) \, \log b -2\log(e/2),
\end{eqnarray*}
where $\,N=K^2L$,
$\,D=[\Q(\alpha_1,\alpha_2,\alpha_3) :
\Q]\bigm/[\R(\alpha_1,\alpha_2,\alpha_3) : \R]$,
\[
g=\frac{1}{4}- \frac{N}{12RST}, \qquad
b = (b_2\eta_0) (b_2\zeta_0)\left(\prod_{k=1}^{K-1}k!\right)^{-\frac{4}{K(K-1)}},
\]
where
\[
\eta_0 = \frac{R-1}{2}+\frac{b_1}{b_2} \times \frac{S-1}{2},\qquad
\zeta_0 = \frac{T-1}{2}+\frac{b_3}{b_2} \times \frac{S-1}{2},
\]
and
\[
a_i \ge  \rho|\log \alpha_i|- \log |\alpha_i|+2D \h(\alpha_i), \qquad
\text{$i=1$, 2, 3.}
\]
{\bf If}
\[
4(R_1+1)(S_1+1) \ge T_1+1,
\leqno{\rm (i)}
\]
\[
4(R_1+1)(T_1+1) \ge S_1+1,
\leqno{\rm (ii)}
\]
\[
 (R_1+1)(S_1+1)(T_1+1) \ge 12(K-1)^2(L-1),
\leqno{\rm (iii)}
\]
and
\[
4(R_1+1)(S_1+1)(T_1+1) \ge 8(2K+L-2)^2
\leqno{\rm (iv)}
\]
 {\bf then either}
\[
 \Lambda' > \rho^{-KL}
\]
 where
\[
\Lambda'=|\Lambda| \cdot
\max \left\{
\frac{LRe^{LR |\Lambda|/(2b_1)}}{2b_1},
\frac{LSe^{LS |\Lambda|/(2b_2)}}{2b_2},
\frac{LTe^{LT |\Lambda|/(2b_3)}}{2b_3}
\right\}
\]
{\bf or}
 at least one of the conditions ({C1}), ({C2}), ({C3}) of Proposition~\ref{prop:Ma4}
hold.
\end{theorem}
\section{Proof of Theorem~\ref{thm:Fibonacci}} 
We are now ready to complete the proof of Theorem~\ref{thm:Fibonacci}.
We argue by contradiction. Suppose that there is a perfect
power in the Fibonacci sequence other than those listed
in Theorem~\ref{thm:Fibonacci}.
By Propositions~\ref{prop:FibPrimeIndex} and~\ref{prop:pbig} 
there  is a solution
$(n,y,p)$ to~(\ref{eqn:FibPrimeIndex}) with $p > 2 \times 10^8$. 

Recall that $F_n=(\omega^n-\omega^{-n})/{\sqrt5}$. Thus the linear
form
\[
\Lambda = n\log \omega - \log \sqrt5 -p\log y
\]
satisfies
\[
\log |\Lambda | < -2p \log y+1.
\]
By Proposition~\ref{prop:Fib733} 
\[
\log y >10^{20}
\]
(and indeed much more).
It seems very difficult to get good lower bounds for $|\Lambda |$ when
it is written in the
previous form. We write
\[
n = k p -q, \qquad \text{where \ $0\le q <p$,}
\]
[notice that $q$ is not necessarily a prime number,
but we have some lack of letters!]. Then
\[
\Lambda = p\log (\omega^k/y) - q\log\omega - \log \sqrt5
\]
and it is easy to see that it is now of the right form.
We know that $p>2 \times 10^8$ and we will obtain a contradiction
by showing that $p < 2 \times 10^8$ using our Theorem~\ref{thm:Maurice}.

The first step is to get an upper bound on $p$ free of any condition.
For this purpose our
Theorem~\ref{thm:Maurice} is inconvenient to use: we have to deal with the conditions (C1),
(C2) and (C3). This is
the reason why we first apply Matveev's estimate (Corollary 2.3): assume
$\Lambda\not=0$,
if real numbers $A_j$ satisfy
\[
A_j \ge \max\bigl\{D \h(\alpha_j), |\log \alpha_j|, 0{.}16\bigr\},\quad
1\le j \le 3,
\]
and if
\[
B=\max\{|b_1|, |b_2|,|b_3|\}
\]
then
\[
\log |\Lambda | > \frac{3e}{2} 30^{6}3^{3{.}5} D^2 A_1A_2A_3 \log(eD)
\log(eB),
\]
where $D=2$ and $B=p$ in our case. This leads to
\[
p<2{.}4\times 10^{13}.
\]

We can now apply Theorem~\ref{thm:Maurice} with
\[
\alpha_1 = \omega , \quad \alpha_2 =\omega^k/y, \quad \alpha_3 =\sqrt5,
\quad D=2
\]
and
\[
b_1 = q , \quad b_2 =p, \quad b_3 =1.
\]
We can take
\[
a_1=(\rho+3)\log \sqrt5, \quad a_2=  (\rho +2p)\log \omega + 4 \log y>4\cdot 10^{20},\quad
a_3=(\rho+1)\log \omega,
\]
where $\rho>e$,
notice that the condition $a_3\le a_1$ is satisfied.
To apply the Theorem, we shall choose some rational integer
\[
L\ge 100
\]
and put
\[
K=\lfloor m L a_1 a_2 a_3 \rfloor , \quad \hbox {with \ $10<m<50$}
\]
and
\[
R_1= \lfloor c_1 L^{2/3}  a_2 a_3 \rfloor, \quad
S_1= \lfloor c_1 L^{2/3} a_1  a_3 \rfloor, \quad
T_1= \lfloor c_1 L^{2/3}  a_1 a_2 \rfloor, \quad \text{with \
$c_1=(32{.}001\,m^2)^{1/3}$}
\]
and
\[
R_2= \lfloor c_2 L   a_2 a_3 \rfloor, \quad
S_2= \lfloor c_2 L  a_1  a_3 \rfloor, \quad
T_2= \lfloor c_2 L   a_1 a_2 \rfloor, \quad \text{with \
$c_2=(12\,m^2)^{1/3}$}
\]
and we put
\[
R=R_1+R_2+1, \quad
S=S_1+S_2+1, \quad
T=T_1+T_2+1. \quad
\]
With such a choice it is easy to check that the four conditions (i),
(ii), (iii), (iv) hold.
And we get
\[
\log |\Lambda|> - KL\log \rho -\log(KL)
\]
or at least one of the conditions ({C1}), ({C2}) and ({C3}) hold.
First notice that, in our case (where $b_2=p$ is prime)
\[
\bigl(({C1})\ {\rm or}\ ({C2})\bigr) \ \Rightarrow \ p \le
\max\{S_1,S_2\}.
\]
Thus, if $p > \max\{S_1,S_2\}$ then ({C1}) and ({C2}) do not hold and
then (C3) holds.
If (C3) holds then recall that
\[
 s't'b_1+r't''b_2+r's'b_3=0
\]
 where the factors of the $b_i$'s are bounded above as in
Proposition~\ref{prop:Ma4}. In the previous
 relation we may assume that $r'$ and $s'$ are coprime, then
 $s'$ divides $t''b_2$. For us, $b_2=p$ and $|s'|<p$, thus $s'$ divides
$t''$ and we get
\[
 t' b_1+r't_2b_2+ r'b_3=0, 
\]
 for some integer $t_2$.
 We may also assume that $r'$ and $t'$ are coprime, which implies
 that $r'$ divides $b_1$, say $b_1=r'b'$, and we have the relation
  --- since here $b_3=1$ ---  we obtain the relation
\[
 t' b'+t_2p+ 1=0, \qquad \text{with \ $q=r'b'$ and $t_2$ divides $t''$}
\]
 where
\[
   |t_2|\le |t'|\le 1+\left(\frac{(S_1+1)(T_1+1)}{R_1+1 } \right)^{1/2}
 < 1+1{.}0001 \cdot \left( \frac{(S_1+1)a_1}{a_3} \right)^{1/2}
\]
 and
\[
 0<r'<\left(\frac{(R_1+1)(S_1+1)}{T_1+1 } \right)^{1/2}
 <1{.}0001 \cdot \left( \frac{(S_1+1)a_3}{a_1} \right)^{1/2}.
\]
In this case, we rewrite $\Lambda$ as a linear form in two logarithms:
\[
\Lambda = p \log \bigl(\alpha_2 {\alpha_3}^{t_2}\bigr)
-b'\log \bigl(\omega^{r'} {\alpha_3}^{t'}\bigr)
\]
where
\[
\omega = \frac{1+\sqrt 5}{2}, \qquad
\alpha_2 = \frac{\omega^k}{y}, \qquad
\alpha_3 = \sqrt 5.
\]
This ends our preliminary discussion.
Now, after some computer search we see that we can apply 
Theorem~\ref{thm:Maurice} with the choices
\[
  L=260, \qquad
  \rho=11, \qquad  m=26{.}12446,
\]
and then, in the first case, we get
\[
p < 451 \times 10^6.
\]
With this choices, we have
\[
\max\{S1,S2\}=\max\{63054,290211\}<10^8
\]
so that neither condition (C1) nor condition (C2) hold. And,
\[
|r'|\le 179, \qquad |t'|,\, |t_2| \le 354.
\]

We apply the following result, Corollaire 2 of \cite{LMN} and Tableau 3:
\begin{lem}\label{lem:Ma5}
Let $\alpha_1$ and $\alpha_2$ be positive real algebraic numbers
which are multiplicatively independent.
Let
\[
\Lambda =b_1\alpha_1-b_2\alpha_2
\]
where $b_1$ and $b_2$ are positive rational integers.
Put
\[
D=[\Q(\alpha_1,\alpha_2):\Q]
\]
and let $A_1$, $A_2$ be real numbers ${}>1$ such that
\[
\log A_i \ge \max\left\{ \h(\alpha_i), \frac{|\log \alpha_i|}{D} ,
\frac{1}{D}\right\},
\qquad i=1,\,2.
\]
Put also
\[
b'=\frac{b_1}{D \log A_2}+ \frac{b_2}{D \log A_1}.
\]
Then
\[
\log |\Lambda | \ge - 25{.}55\,D^4\,
\left( \max\left\{ \log b'+0{.}19,\, \frac{18}{D}, 1\right\}\right)^2 \,
\log A_1\,\log A_2.
\]
\end{lem}
Here, with our initial notation ({\it i.e.} $\alpha_2=\omega^k/y$,
$\alpha_3=\sqrt5$),
 we have $D=2$ and choose (as we may!)
\[
\log A_1 = 1{.}001\,\h(\alpha_2)
\]
so that
\[
\log A_1\ge \h(\alpha_2)+|t_2 |\h(\alpha_3)+1 >10^{20}
\]
and
\[
  \log A_2 = 354\,\h(\alpha_3)+179\, \h(\omega)+1=328{.}93896\dots
\]
Then
\[
b'=\frac{p}{D \log A_2}+ \frac{b'}{D \log A_1}< \frac{p}{2}
\left( \frac{1}{10^{20}}+ \frac{1}{\log A_2 }\right)< \frac{p}{657}. 
\]
and
\[
\log |\Lambda | \ge - 25{.}55\,D^4\,
\bigl( \max \bigl\{ \log (p/657)+0{.}19,\, 9\bigr\}\bigr)^2
\, \log A_1\,\log A_2.
\]
and we get (recall that we know $p>2 \times 10^8$)
\[
\log |\Lambda | \ge - 25{.}6\times 16 \times 329 \times
\bigl( \max \bigl\{ \log (p/657)+0{.}19,\, 9\bigr\}\bigr)^2\times \log
y
\]
and, when compared to the upper bound of $\log |\Lambda |$, this leads
to
\[
p < 7\times 10^7.
\]
\medskip
Thus we have proved that $p<451\times 10^6$. From this upper bound, if
we iterate five times
this process, we get:
\[
p < 183 \times 10^6.
\]
which certainly shows that $p< 2 \times 10^8$ and we have obtained our
contradiction. This completes the proof of Theorem~\ref{thm:Fibonacci}.


\begin{thebibliography}{}
\bibitem{BW} A. Baker and W\"{u}sholz,
{\em Logarithmic forms and group varieties},
J. reine angew. Math. {\bf 442} (1993), 19--62.


\bibitem{PARI} C. {Batut}, K. {Belabas}, D. {Bernardi}, H. {Cohen} and M. {Olivier}, 
{\em User's guide to PARI-GP}, 
version 2.1.1. (See also {\tt http://www.parigp-home.de/}.)

\bibitem{B} C.D. Bennett, J. Blass, A.M.W. Glass, D.B. Meronk, R.P. Steiner,
{\em Linear forms in the logarithms of three positive rational numbers},
J. Th\'eor. Nombres Bordeaux {\bf 9} (1997), 97--136.

\bibitem{Be} M. A. Bennett, 
{\em Rational approximation to algebraic number of
small height: The Diophantine equation} $\mid a x^n - by^n \mid =1$, 
J. reine angew. Math. 535 (2001), 1--49.

\bibitem{BS} M. A. Bennett and C. M. Skinner, 
{\em Ternary Diophantine equations via Galois representations and modular forms},
to appear in the Canadian Journal of Mathematics.

\bibitem{BHV} Yu. Bilu, G. Hanrot and P. M. Voutier,
{\em Existence of primitive divisors of Lucas and Lehmer numbers}.
With an appendix by M. Mignotte. J. reine angew. Math.
{\bf 539} (2001), 75--122.

\bibitem{Magma} W.~Bosma, J.~Cannon and C.~Playoust: {\em The Magma
Algebra System I: The User Language}, J. Symb. Comp. {\bf 24} (1997),
235--265.  (See also {\tt http://www.maths.usyd.edu.au:8000/u/magma/}.)

\bibitem{Mod} C. Breuil, B. Conrad, F. Diamond, R. Taylor, 
{\em On the modularity of elliptic curves over $\Q$: wild $3$-adic exercises},
J. Amer. Math. Soc. {\bf 14 No.4} (2001), 843--939.


\bibitem{BGunits} Y. Bugeaud and K. Gy\H ory,
{\em Bounds for the solutions of unit equations},
Acta Arith. {\bf 74} (1996), 67--80.

\bibitem{BG} Y. Bugeaud and K. Gy\H ory, 
{\em Bounds for the solutions of Thue-Mahler equations and norm form equations}, 
Acta. Arith. {\bf 74} (1996), 273--292.

\bibitem{BM} Y. Bugeaud and M. Mignotte,  
{\em On integers with identical digits},
Mathematika {\bf 46} (1999), 411--417.

\bibitem{BMRS} Y. Bugeaud, M. Mignotte, Y. Roy and T. N. Shorey, 
{\em The equation $(x^n-1)/(x-1)=y^q$ has no solutions with $x$ square}, 
Math. Proc. Cambridge Phil. Soc. {\bf 127} (1999), 353--372. 

\bibitem{BMS} Y. Bugeaud, M. Mignotte and S. Siksek,
{\em Classical and modular approaches to exponential Diophantine
equations. II}, in preparation.


\bibitem{Cohn1} J. H. E. Cohn, 
{\em On square Fibonacci numbers},
J. London Math. Soc. {\bf 39} (1964), 537--540.

\bibitem{Cohn3} J. H. E. Cohn,
{\em Lucas and Fibonacci numbers and some Diophantine equations},
Proc. Glasgow Math. Assoc. {\bf 7} (1965), 24--28.

\bibitem{Cohn2} J. H. E. Cohn, 
{\em Perfect Pell powers},
Glasgow Math. J. {\bf 38} (1996), 19--20.

\bibitem{Cre} J. E. Cremona,
{\em Algorithms for modular elliptic curves},
2nd edition, Cambridge University Press, 1996.


\bibitem{DM} H. Darmon and L. Merel,
{\em Winding quotients and some variants of Fermat's
Last Theorem}
J. reine angew. Math. {\bf 490} (1997), 81--100.



\bibitem{Iv} W. Ivorra,
{\it Courbes elliptiques sur $\Q$, ayant un point d'ordre $2$ 
rationnel sur $\Q$, de conducteur $2^N p$}. Preprint.


\bibitem{Kra} A. Kraus, 
{\em Sur l'\'equation $a^3+b^3=c^p$},
Experimental Mathematics {\bf 7} (1998), No. 1, 1--13.

\bibitem{KO} A. Kraus and J. Oesterl\'e, 
{\em Sur une question de B. Mazur}, 
Math. Ann. {\bf 293} (1992), 259--275.

\bibitem{Landau} E. Landau, {\em Verallgemeinerung eines
P\'olyaschen Satzes auf algebraische Zahlk\"orper},
Nachr. Kgl. Ges. Wiss. G\"ottingen, Math.-Phys. Kl. (1918), 478--488.


\bibitem{La} M. Laurent, 
{\em Linear form in two logarithms and the interpolation determinants}, 
Acta Arith. {\bf 66} (1994), 181--199.

\bibitem{LMN} M. Laurent, M. Mignotte and Y. Nesterenko,
{\em Formes lin\'eares en deux logarithmes et d\'eterminants d'interpolation}, 
J. Number Theory {\bf 55} (1995), 255--265.

\bibitem{Len} H. W. Lenstra, Jr.,
{\em Algorithms in algebraic number theory},
Bull. Amer. Math. Soc. {\bf 26} (1992), 211-244.





\bibitem{LF} H. London and R. Finkelstein, 
{\em On Fibonacci and Lucas numbers which are perfect powers}, 
Fibonacci Quart. {\bf 5} (1969), 476--481.

\bibitem{Mc} J. McLaughlin,
{\em Small prime powers in the Fibonacci sequence},
to appear.

\bibitem{Mat} Matveev, {\em An explicit lower bound for a
homogeneous rational linear form in logarithms of algebraic numbers. II},
Izv. Ross. Akad. Nauk Ser. Mat. {\bf 64} (2000), 125--180.
English transl. in Izv. Math. {\bf 64} (2000), 1217--1269.



\bibitem{Mazur} B. Mazur, {\em Rational isogenies of prime degree},
Invent. Math. {\bf 44} (1978), 129--162.

\bibitem{Mi3} M. Mignotte, 
{\em Entiers alg\'ebriques dont les conjugu\'es sont 
proches du cercle unit\'e}, 
S\'eminaire Delange--Pisot--Poitou, 19e ann\'ee: 1977/78, 
Th\'eorie des nombres, Fasc. 2, Exp. No. 39, 6 pp., Secr\'etariat 
Math., Paris, 1978. 

\bibitem{Nark} W. Narkiewicz,
{\em Elementary and Analytic Theory of Algebraic Numbers,} 
Springer-Verlag, Berlin, 1990.

\bibitem{Pe82} A. Peth\H{o},
{\em Perfect powers in second order linear recurrences},
J. Number Theory {\bf 15} (1982), no. 1, 5--13.

\bibitem{Pe83} A. Peth\H{o},
{\em Full cubes in the Fibonacci sequence},
Publ. Math. Debrecen {\bf 30} (1983), no. 1, 117--127.

\bibitem{Pe01} A. Peth\H{o},
{\em Diophantine properties of linear recursive sequences II},
Acta Math. Paedagogicae Ny\'{i}regyh\'{a}ziensis {\bf 17} (2001), 81--96.


\bibitem{Ribet} K. Ribet, {\em On modular representations of
${\rm Gal}(\overline{\Q}/\Q)$ arising from modular forms}, Invent. Math.
{\bf 100} (1990), 431--476.

\bibitem{Rob} N. Robbins,
{\em On Fibonacci numbers which are powers. II}.
Fibonacci Quart. {\bf 21} (1983), no. 3, 215--218.

\bibitem{ShSt} T. N. Shorey and C. L. Stewart,
{\em On the Diophantine equation $a x^{2t} + b x^t y + c y^2 = d$ and
pure powers in second order linear recurrences},
Math. Scand. {\bf 52} (1983), 24--36.


\bibitem{ST} T. N. Shorey and R. Tijdeman,
{\em Exponential Diophantine equations}, 
Cambridge Tracts in Mathematics 87,
Cambridge University Press, Cambridge, 1986.

\bibitem{Siegel} C. L. Siegel, 
{\em Absch\"atzung von Einheiten}, Nachr. Akad. Wiss. G\"ottingen II,
Math.-Phys. Kl., Nr. 9, (1969), 71--86.

\bibitem{SC} S. Siksek and J. E. Cremona, 
{\em On the Diophantine equation $x^2+7=y^m$}, 
Acta Arith. {\bf 109.2} (2003), 143--149.


\bibitem{TW} R.L. Taylor and A. Wiles, 
{\em Ring theoretic properties of certain Hecke algebras},
Annals of Math. {\bf 141} (1995), 553--572.

\bibitem{Vou} P. M. Voutier, {\em An effective lower bound for the
height of algebraic numbers}, Acta Arith. {\bf 74} (1996), 81--95.


\bibitem{Wa} M. Waldschmidt, 
{\em Minorations de combinaisons lin\'eaires de logarithmes de nomb\-res alg\'e\-bri\-ques}, 
Canadian J. Math.  {\bf 45} (1993), 176-224.

\bibitem{W} M. Waldschmidt, 
{\em Diophantine approximation on linear algebraic groups}, 
Springer, Berlin, 2000.

\bibitem{Wi} A. Wiles, 
{\em Modular elliptic curves and Fermat's Last Theorem},
Annals of Math. {\bf 141} (1995), 443--551.

\bibitem{Wy} O. Wyler, 
{\em Solution to Problem 5080}, 
Amer. Math. Monthly {\bf 71} (1964), 220--222.
\end{thebibliography}
\end{document}